\documentclass[12pt]{amsart}
\usepackage{amsmath,amsfonts,amssymb,amscd,amsthm,amsbsy,epsf}
\newtheorem{thm}{Theorem}[section]
\newtheorem*{thmA}{Theorem A}
\newtheorem*{thmB}{Theorem B}
\newtheorem{lem}[thm]{Lemma}
\newtheorem{cor}[thm]{Corollary}
\newtheorem{prop}[thm]{Proposition}
\theoremstyle{definition} 		%% bold "label" with roman text
\newtheorem{defn}[thm]{Definition}
\newtheorem{remark}[thm]{Remark}
\def\nat{{\mathbb N}}
\def\A{{\mathcal A}}
\def\B{{\mathcal B}}
\def\C{{\mathcal C}}
\def\D{{\mathcal D}}
\def\F{{\mathcal F}}
\def\G{{\mathcal G}}
\def\P{{\mathcal P}}
\def\U{{\mathcal U}}
\def\bs{{\mathbf s}}
\def\bt{{\mathbf t}}

\def\Tau{{\mathfrak T}}
\begin{document}
\title{Block Combinatorics}
\author{V. Farmaki and S. Negrepontis}
%\dedicatory{}
\begin{abstract}
In this paper we extend the block combinatorics partition theorems of 
Hindman and Milliken in the setting of the recursive system of the block
Schreier families $(\B^\xi)$ consisting of families defined for every 
countable ordinal $\xi$. 
Results contain 
(a)~a block partition Ramsey theorem for every countable 
ordinal $\xi$ (Hindman's theorem corresponding to $\xi=1$, and Milliken's 
theorem to $\xi$ a finite ordinal), 
(b)~a countable ordinal form of the block Nash-Williams partition theorem, 
and (c)~a countable ordinal block partition theorem for sets closed in the 
infinite block analogue of Ellentuck's topology.
\end{abstract}
%\address{Department of Mathematics, Athens University, Athens 157 84, Greece}
%\address{Department of Mathematics, The University of Texas at Austin,
%Austin, TX 78712 U.S.A.}
%\email{vfarmaki@math.uoa.gr, vfarmaki@math.utexas.edu}
%\email{snegrep@math.uoa.gr, snegrep@math.utexas.edu}
\subjclass{}  
%\thanks{Here is a fake footnote.}
\keywords{block Ramsey, Nash-Williams combinatorics, Schreier families}
\subjclass{Primary 05D10; 46B20}

\maketitle
\baselineskip=18pt
\pagestyle{plain}              % page nos. at bottom, no headline

%%%%%%   BODY   %%%%%%%%%%%%%%%%%%%%%%%%%%%%%%

\section*{Introduction}

In this paper we extend the block combinatorics partition theorems of 
Hindman \cite{H} and Milliken \cite{M} in the setting of the recursive 
system of the block Schreier families $(\B^\xi)_{\xi<\omega_1}$, 
consisting of families defined for every countable ordinal $\xi$. 
The main results are contained in Theorem~1.8, Theorem~3.3 and 
Corollary~3.6, and Theorem~4.2. 
Theorem~1.8 is a block Ramsey partition theorem for every countable 
ordinal $\xi$. 
Hindman's theorem (\cite{H}) corresponds to $\xi =1$, and Milliken's 
theorem (Theorem~2.2 in \cite{M}) corresponds to $\xi$ a finite ordinal 
(natural number). 
Using this result, together with  some additional properties of the Schreier 
system, including its Cantor-Bendixson index, we then obtain strong 
countable ordinal forms of the block Nash-Williams partition theorem 
(Theorem~3.3 and Corollary~6.6). 
A special case of these results is Milliken's theorem (Theorem~3.5 in 
\cite{M}). 
A consequence of Theorem~3.3 gives Theorem~4.2, a countable ordinal block 
partition for sets closed in the infinite block analogue of Ellentuck's 
topology. 
A special case of this result is Milliken's Theorem~4.4 in \cite{M}. 

The recursive system of block Schreier families $(\B^\xi)_{\xi<\omega_1}$, 
each family consisting of finite blocks of finite subsets of $\nat$, introduced 
in the present paper (Definition~1.3), constitutes the block analogue of the 
recursive system of Schreier families $(\A_\xi)_{\xi<\omega_1}$, 
each family consisting of finite subsets of $\nat$, defined and studied in 
\cite{F1}, \cite{F2}, \cite{F3}. 
Although the recursive system of Schreier families $(\A_\xi)$ is a purely 
combinatorial entity, intimately related, as it turns out, with Ramsey 
and Nash-Williams combinatorics, it nevertheless arose gradually in 
connection with the theory of Banach  spaces, originally by Schreier 
(\cite{S}) (for $\xi=\omega$), next by Alspach-Odell \cite{AO} (for 
$\xi =\omega^\kappa$, $\kappa$ a natural number) and Alspach-Argyros
\cite{AA} (for $\xi=\omega^\alpha$, $\alpha$ a countable ordinal),
and finally by 
Farmaki \cite{F1}, \cite{F2}, \cite{F3} and Tomczak-Jaegermann \cite{TJ} 
(for $\xi$ any countable ordinal).

The combinatorial nature of the block Schreier families 
$(\B^\xi)_{\xi<\omega_1}$ and of the block combinatorics developed in this 
paper makes it reasonable to expect that the results presented in this 
paper will find substantial applications in various branches of mathematics.
\bigskip

We now turn into a somewhat more detailed description of the results of 
this paper. 
The statement of the general block partition Theorem~1.8 is the following:

\begin{thmA}
Let $\F$ be a subfamily of th family $\B^{<\omega}(\nat)$ of all disjoint 
collections on $\nat$, $\D_0$ an infinite disjoint collection on $\nat$ 
and $\xi$ a countable ordinal.
Then there exists an infinite disjoint collection $\D <\D_0$ such that 
\begin{itemize}
\item[{}] either $\B^\xi \cap \B^{<\omega} (\D) \subseteq \F$ or 
$\B^\xi \cap \B^{<\omega}(\D)\subseteq\B^{<\omega} (\nat)\setminus \F$.
\end{itemize}
\end{thmA}

As already mentioned, the case $\xi=1$ Theorem~A is Hindman's theorem, 
and the case $\xi$ is a natural number is Milliken's theorem. 
Theorem~A, holding for every countable ordinal, in terms of the block 
Schreier system, generalizes Hindman's and Milliken's theorems, in an 
analogous way that Ramsey's partition theorem, holding for every countable 
ordinal, in terms of the Schreier system, proved in \cite[Theorem~A]{F2}, 
\cite[Theorem~1.5]{F3}, generalizes the classical Ramsey's theorem 
for $k$-tuples.

The general block Ramsey partition theorem implies strong forms, in terms 
of the block Schreier families $(\B^\xi)_{\xi<\omega_1}$, of block 
Nash-Williams partition theorem: 
Theorem~3.1 for a partition family $\F$ that is hereditary, 
Theorem~3.3 for $\F$ a tree, and Corollaries~3.4, 3.6. 
The tools for proving these results include the thinness of the  recursive 
Schreier families (Proposition~2.2), the canonical representation of every, 
finite or infinite, disjoint collection with respect to every family 
$\B^\xi$ (Proposition~2.4), and the concept of the strong block 
Cantor-Bendixson index and its computation for the family $\B^\xi$ 
(Proposition~2.10). 

The statement of Theorem 3.3, in a somewhat abbreviated form, is the following:

\begin{thmB}
Let $\F$ be a subfamily of the family $\B^{<\omega}(\nat)$ of all disjoint 
collections on $\nat$, which is a tree, and let $\D$ be an infinite disjoint 
collection on $\nat$ and denote by $\F_h$ the maximal hereditary subfamily 
of $\F$ contained in $\F$. 
Then we have the following cases: 

\noindent {\bf [Case 1]}
The family $\F_h\cap \B^{<\omega} (\D)$ is not pointwise closed.

Then, there exists an infinite disjoint collection $\D_1<\D$ such that
$\B^{<\omega} (\D_1) \subseteq \F$.

\noindent {\bf [Case 2]}
The family $\F_h \cap \B^{<\omega} (\D)$ is pointwise closed.
Setting
\begin{equation*}
\zeta_{\D}^{\F} =  \xi_{\D}^{\F_h} = \sup
\{ s_{\D_1} (\F_h) : \D_1 <\D\}\ ,
\end{equation*}
which is a countable ordinal, the following subcases obtain:
\begin{itemize}
\item[2(i)]
If $\xi+1 <\zeta_{\D}^{\F}$, then there exists $\D_2<\D$
such that
$$\B^\xi \cap \B^{<\omega} (\D_2) \subseteq \F\ ;$$
\item[2(ii)]
if $\xi+1>\xi>\zeta_{\D}^{\F}$, then for every $\D_1<\D$
there exists $\D_2<\D_1$  such that
$$\B^\xi\cap \B^{<\omega} (\D_2) \subseteq \B^{<\omega} (\nat)\setminus \F$$
(equivalently $\F\cap \B^{<\omega} (\D_2) \subseteq(\B^\xi)^*\setminus
\B^\xi$). 
%\item[2(iii)]
%if $\xi+1 = \zeta_{\D}^{\F}$ or $\xi = \zeta_{\D}^{\F}$, then there
%exists $\D_2 <\D$ such that either
%${\B_\xi\cap \B^{<\omega} (\D_2)\subseteq \F}$ or
%$\B_\xi\cap \B^{<\omega} (\D_2)\subseteq \B^{<\omega}(\nat)\setminus\F$.
\end{itemize}
\end{thmB}

Generally the results of Section 3 are more detailed versions of the 
block Ramsey Theorem~1.8, made possible by additional conditions for 
the partition $\F$. 
Roughly speaking the given partition $\F$ (assumed to be hereditary, or a 
tree), is strong enough to absorb and contain all $\B^\xi$ finite blocks 
of a suitably defined infinite disjoint collection, in case the strong 
Cantor-Bendixson index of $\F$ exceeds the ordinal $\xi$, and is unable 
to absorb these $\B^\xi$ blocks and obliged to allow its complement 
$\F^c$ to receive a substantial portion of them (albeit, in general in an 
asymmetrically weaker manner), in case the index of $\F$ falls below the 
ordinal $\xi$.

The reason for regarding the block Ramsey partition theorems for trees 
of Section~3 as strong block Nash-Williams partition theorems may not be 
immediately apparent to the reader, but it is clearly supported by the fact 
that these theorems produce effortlessly results (such as Corollary~3.8) 
of the more traditional Nash-Williams form.

Finally we prove, in Theorem~4.2, a strong countable ordinal version 
involving the block Schreier families $\B^\xi$, of the block partition 
theorem for sets closed in the block  analogue of the Ellentuck topology 
on all infinite blocks. 
Theorem~4.2  contains as a special case (Corollary~4.7) Milliken's 
Theorem~4.4 in \cite{M}, an analogue for blocks of Ellentuck's 
theorem (\cite{E}).

We will need the following notations.

\subsection*{Notation}
We denote by $\nat$ the set $\{1,2,\ldots,\}$ of all natural numbers. 
For an infinite subset $M$ of $\nat$, we denote by $[M]^{<\omega}$ the set 
of all the finite subsets of $M$ and by $[M]$ the set of all the infinite 
subsets of $M$ (considered as strictly increasing sequences). 

If $s,t$ are non-empty subsets of $\nat$ and $s$ is finite we write $s<t$ 
if $\max s <\min t$, while $s\propto t$ if $s$ is an initial segment of $t$. 

A collection (finite or infinite) of disjoint finite subsets of $\nat$ is 
called a disjoint collection. 
We denote by $\D_1,\D_2,\ldots$ the infinite disjoint collections (considered 
as strictly increasing sequences of finite subsets of $\nat$). 
For an infinite disjoint collection $\D$ (resp.\ for a finite disjoint 
collection $\bs$) we denote by $FU(\D)$ (resp.\ by $FU(\bs)$) the family 
of all the finite unions of elements of $\D$ (resp.\ of $\bs$). 
For two infinite disjoint collections $\D_1,\D_2$ we write $\D_1<\D_2$ 
if $\D_1\subseteq FU(\D_2)$. 
If $\bs$ is a finite collection we write 
$\bs \propto \D_1$ if $\bs$ is an initial segment of $\D_1$ and analogously 
$\bs \propto \bt$, for a finite disjoint collection $\bt$, if $\bs$ is an 
initial segment of $\bt$.

We denote by $\B^{<\omega}(\nat)$ the set of all finite disjoint collections 
and by $\B^\omega (\nat)$ the set of all infinite disjoint collections. 
Thus 
\begin{itemize}
\item[{}] $\B^{<\omega}(\nat) = \{\bs = (s_1,\ldots,s_k) :k\in\nat$, 
$s_1<\cdots < s_k\in [\nat]^{<\omega} \} \cup \{\emptyset\}$;
\item[{}] $\B^\omega (\nat) = \{\D = (s_n)_{n\in\nat} :s_n\in [\nat]^{<\omega}$
and  $s_n< s_{n+1} \ \forall\ n\in\nat\}$.
\end{itemize}
For $\D\in \B^\omega (\nat)$ we denote analogously 
\begin{itemize}
\item[{}] $\B^{<\omega} (\D) = \{\bs = (s_1,\ldots, s_k):k\in\nat$, 
$s_1<\cdots < s_k\in FU(\D)\}\cup \{\emptyset\}$; 
\item[{}] $\B^\omega (\D) = \{\D_1 \in \B^\omega (\nat) : \D_1<\D\}$.
\end{itemize} 

For a family $\A$ of finite subsets of $\nat$, $n\in\nat$ and 
$t\in [\nat]^{<\omega}$, we set 
\begin{itemize}
\item[{}] $\A(n) = \{s\in [\nat]^{<\omega} : \{n\} < s$ 
and $\{n\}\cup s\in \A\}$;
\item[{}] $\A - t = \{s\in [\nat]^{<\omega}  : t< s\}$;
\item[{}] $\A (t) = \{s\in [\nat]^{<\omega} : t\cup s\in \A\}$.
\end{itemize}

For a family $\F$ of finite disjoint collections and $t\in [\nat]^{<\omega}$ 
we set 
\begin{itemize}
\item[{}] $\F(t) = \{\bs = (s_1,\ldots,s_k) \in \B^{<\omega} (\nat) : 
(t,s_1,\ldots, s_k) \in \F\}$;
\item[{}] $\F -t = \{\bs = (s_1,\ldots,s_k) \in \F :t<s_1<\cdots < s_k\}$.
\end{itemize}

%%%%%%%%%%%%%%%%%%%%%
\setcounter{section}{0}
\section{The block Ramsey partition theorem for every countable 
ordinal}

The purpose of this section is to prove the block Ramsey partition  theorem
for every countable ordinal (Theorem~\ref{thm:block-Ramsey}). 
It is the extension to every countable ordinal $\xi$ of Hindman's theorem 
(for $\xi=1$). 

We recall from \cite[Definition 1.3]{F3} the definition of the recursive 
thin Schreier system $(\A_\xi)_{\xi<\omega_1}$ of families of finite 
subsets of $\nat$. 
(In \cite{F3} it was called ``complete'' thin Schreier system.) 
The definition of the system $(\A_\xi)_{\xi<\omega_1}$ employs the 
canonical representation of ordinals (cf.\ \cite{K}). 

\begin{defn}[The recursive thin Schreier system, 
{[F1, Def. 7], [F2, Def. 1.5] [F3, Def. 1.3]}]
%\cite[Definition 1.3]{F}
\label{recursivethin}
For every non-empty limit ordinal $\lambda$ choose and fix $\alpha$ strictly 
increasing sequence $(\lambda_n)$ of successor ordinals smaller than $\lambda$ 
with $\sup \lambda_n = \lambda$. 
The system $(\A_\xi)_{\xi<\omega_1}$ is defined recursively as follows: 
\begin{itemize}
\item[(1)] $\A_0 = \{\emptyset\}$;
\item[(2)] $\A_{\zeta+1} = \{s\in [\nat]^{<\omega} : s= \{n\} \cup s_1$, 
where $n\in \nat$, $\{n\} <s_1$ and $s_1\in \A_\zeta\}$; 
\item[(3)] $\A_{\omega^{\beta+1}} = \{s\in [\nat]^{<\omega} : 
s = \bigcup_{i=1}^n s_i$, where $n= \min s_1$, $s_1<\cdots < s_n$ and 
$s_1,\ldots, s_n\in \A_{\omega^\beta}\}$; 
\item[(4)] for a non-zero, countable limit ordinal $\lambda$, 
\newline
$\A_{\omega^\lambda} = \{s\in [\nat]^{<\omega} : s\in \A_{\omega^{\lambda_n}}$
with $n= \min s\}$; and 
\item[(5)] for a limit ordinal $\xi$ such that $\omega^a< \xi < 
\omega^{a+1}$ for some $0< a<\omega_1$, if \newline
$\xi = p\omega^ a
+ \sum_{i=1}^m p_i \omega^{a_i}$, where $m\in \nat$ with $m\ge0$, 
$p,p_1,\ldots,p_m$ are natural numbers with $p,p_1,\ldots,p_m\ge1$ 
(so that either $p>1$, or $p=1$ and $m\ge 1$) and 
$a,a_1,\ldots,a_m$ are ordinals with $a>a_1>\cdots a_m >0$, 
\newline
$\A_\xi = \{s\in [\nat]^{<\omega} :s= s_0 \cup (\bigcup_{i=1}^m s_i)$ 
with $s_m < \cdots < s_1 <s_0$,
$s_0= s_1^0\cup\cdots\cup s_p^0$ with $s_1^0<\cdots < s_p^0\in \A_{\omega^a}$,
and $s_i = s_1^i \cup\cdots\cup s_{p_i}^i$ with 
$s_1^i <\cdots<\ s_{p_i}^i\in \A_{\omega^{a_i}}$ $\forall\ 1\le i\le m\}$.
\end{itemize}
\end{defn}

We will also need the following results proved in \cite{F3}:

\begin{prop}\label{prop:recursive}
{\rm (i) (Recursiveness of $(\A_\xi)_{\xi<\omega_1}$)} 
For every countable ordinal $\xi$ there exists a concrete sequence 
$(\xi_n)$ of countable ordinals with  $\xi_n <\xi$ such that 
\begin{equation*}
\A_\xi (n) = \A_{\xi_n} \cap [\{n+1,n+2,\ldots\}]^{<\omega}\ \text{ for every }
\ n\in \nat\ .
\end{equation*}
Moreover, $\xi_n = \zeta$ for every $n\in\nat$ if $\xi = \zeta+1$ and 
$(\xi_n)$ is a strictly increasing sequence with $\sup_n \xi_n=\xi$ if $\xi$ 
is a limit ordinal 
{\rm (\cite[Proposition 1.7]{F3}}.

{\rm (ii) (Thinness of $\A_\xi$)}
The families $\A_\xi$, $\xi<\omega_1$, are thin (i.e., there are no elements 
$s,t\in \A_\xi$ with $s\propto t$ and $s\ne t$) 
{\rm (\cite[Proposition 2.2]{F3})}. 

{\rm (iii) (Canonical representation w.r.t.\ $\A_\xi$)}
Every (finite or infinite) non-empty subset of $\nat$ has canonical 
representation with respect to $\A_\xi$ (i.e., if $I\in[\nat]$, then there 
exists a unique sequence $(s_n)_{n\in\nat}\subseteq \A_\xi$ such that 
$I = \bigcup_{n=1}^\infty s_n$ and $s_1<s_2<\cdots$; and 
if $s\in [\nat]^{<\omega}$, then there exist unique $n\in \nat$, 
$s_1,\ldots, s_n\in \A_\xi$ and a, possible empty, set $s_{n+1}$ which is a 
proper initial segment of some element of $\A_\xi$ such that 
$s = \bigcup_{i=1}^{n+1} s_i$) 
{\rm (\cite[Proposition 2.4]{F3})}. 

{\rm (iv) ($\xi$-Ramsey partition theorem)} 
For an arbitrary family $\F$ of finite subsets of $\nat$, $M\in[\nat]$ 
and $\xi<\omega_1$ there exists $L\in [M]$ such that either 
$\A_\xi\cap [L]^{<\omega} \subseteq \F$ or \newline
$\A_\xi\cap [L]^{<\omega} \subseteq [\nat]^{<\omega} \setminus \F$ 
{\rm (\cite[Theorem 1.5]{F3})}. 
\end{prop}

In most of this paper we will be working with finite disjoint collections 
(= finite sets of blocks) instead of finite sets; and, thus we will be 
working with the recursive thin block Schreier system $(\B^\xi)_{\xi<\omega_1}$
(in place of $(\A_\xi)_{\xi<\omega_1}$), defined below. 

\begin{defn} 
[Recursive thin block Schreier system $(\B^\xi)_{\xi<\omega_1}$] 
\label{recursivethinblock}
We define 
$$\B^0 = \{\emptyset\}\ ;$$ 
and for every countable ordinal $\xi\ge1$ 
\begin{equation*}
\B^\xi = \{\bs = (s_1,\ldots,s_k) : k\in\nat,\ 
s_1<\cdots < s_k \in [\nat]^{<\omega}\ \text{ and }\ 
(\min s_1,\ldots, \min s_k) \in \A_\xi\}\ .
\end{equation*}
\end{defn}

\begin{remark}\label{rem1.4}
(i) For $k\in \nat$,
$\B^k = \{(s_1,\ldots,s_k) :s_1<\cdots < s_k\in [\nat]^{<\omega}\}$.  

(ii) $\B^\omega = \{(s_1,\ldots, s_k) : k\in \nat$, 
$s_1 <\cdots< s_k \in [\nat]^{<\omega}$ and $\min s_1 = k\}$.

(iii) $\B^\xi \in \B^{<\omega} (\nat)$ and
$\emptyset \notin \B^\xi$ for $\xi\ge1$.
\end{remark} 

The following proposition justifies the term ``recursive'' in our definition 
of the system $(\B^\xi)_{\xi<\omega_1}$.

\begin{prop}\label{justification}
For every countable ordinal $\xi\ge 1$ there exists a concrete sequence 
$(\xi_n)$ of countable ordinals with $\xi_n<\xi$ such that for every 
$s\in [\nat]^{<\omega}$ with $\min s=n$ 
\begin{equation*}
\B^\xi (s) = \B^{\xi_n}\cap (\B^{<\omega} (\nat) - s)\ .
\end{equation*}
Moreover, $\xi_n =\zeta$ for every $n\in\nat$ if $\xi = \zeta+1$ and 
$(\xi_n)$ is a strictly increasing sequence with $\sup_n \xi_n=\xi$ if $\xi$ 
is a limit ordinal.
\end{prop} 

\begin{proof}
For every countable ordinal $\xi\ge1$, let $(\xi_n)$ be the corresponding 
sequence  defined in Proposition~1.7 in \cite{F3}
(see Proposition~\ref{prop:recursive}~(i)).

For $\xi=1$, we have $\B^1  = \{(s_1) :s_1\in [\nat]^{<\omega}\}$, hence, 
$\B^1 (s) = \B^0 = \{\emptyset\}$ for every $s\in [\nat]^{<\omega}$. 
Let $1<\xi<\omega_1$ and $s\in [\nat]^{<\omega}$ with $\min s=n$. 
If $(s_1,\ldots,s_k)\in \B^\xi(s)$, then $(s,s_1,\ldots,s_k)\in \B^\xi$, 
hence $(n,\min s_1,\ldots, \min s_k)\in \A_\xi$. 
Thus, by Proposition~1.2 (i), $(\min s_1,\ldots,\min s_k)\in \A_{\xi_n}$ 
and consequently $(s_1,\ldots,s_k)\in \B^{\xi_n}\cap (\B^{<\omega}
(\nat) - s)$. 
On the other hand, if $(s_1,\ldots,s_k)\in \B^{\xi_n}\cap (\B^{<\omega}
(\nat) -s)$, then $(\min s_1,\ldots,\min s_k) \in \A_{\xi_n} \cap 
{[\{n+1,n+2,\ldots \}]^{<\omega} \subseteq \A_\xi (n)}$. 
Hence, $(n,\min s_1,\ldots,\min s_k)\in \A_\xi$ and then 
$(s,s_1,\ldots,s_k)\in \B^\xi$. 
\end{proof}

In the classical  Ramsey theory  (developed to its full extent in \cite{F3}),
the basic, starting dichotomy, corresponding to ordinal level $\xi=1$, is 
the statement that if we partition in two parts an infinite subset of 
$\nat$, then one part at least is infinite. 
For the block Ramsey theory, that we are about to develop, the basic starting 
dichotomy, corresponding to ordinal level $\xi=1$, is 
(the highly non-trivial) Hindman's theorem (\cite{H}) 
(proved also by Baumgartner (\cite{B}), which we now state:

\begin{thm}
[Hindman's theorem, \cite{H}, \cite{B}] 
\label{thm:hindman}
Let $\F \subseteq \B^{<\omega} (\nat)$ be a family of finite disjoint 
collections and $\D_0$ an infinite disjoint collection. 
Then, there exists an infinite disjoint collection $\D<\D_0$ such that either 
$\B^1 \cap \B^{<\omega} (\D)\subseteq \F$ or 
${\B^1 \cap \B^{<\omega} (\D)\subseteq \B^{<\omega}(\nat)\setminus \F}$.
\end{thm}

For the proof of the block-Ramsey partition theorem we will make 
use of a diagonal argument, contained in the following lemma.

\begin{lem}\label{lem:block-Ramsey}
Let $\D_0$ be an infinite disjoint collection and \newline
$\Pi = \{(s,\D):s\in [\nat]^{<\omega}$, $\D\in \B^\omega (\nat)$ with 
$\D<\D_0$ and $s<\bigcup_{F\in \D} F\}$.\newline
If a subset $\G$ of $\Pi$ satisfies:
\begin{itemize}
\item[(i)] for every $(s,\D)\in\Pi$ there exists $(s,\D_1)\in\G$ with 
$\D_1 <\D$; and
\item[(ii)] for every $(s,\D_1) \in\G$ and $\D_2<\D_1$ we have $(s,\D_2)\in\G$,
\end{itemize}
then there exists an infinite disjoint collection $\D<\D_0$ such that 
$(s,\D')\in \G$ for every $s\in FU(\D)$ and $\D' <\D-s$.
\end{lem}

\begin{proof} 
Let $s_0\in \D_0$. 
According to conditions (i) and (ii) there exists $\D_1<\D_0-s_0$ such that 
$(s_0,\D_1)\in\G$ and $(\emptyset,\D_1) \in\G$. 
Let $s_1\in\D_1$. 
Of course, $s_0<s_1$. 
We assume that there have been constructed infinite disjoint collections 
$\D_0,\D_1,\ldots,\D_n$ with $\D_n<\cdots <\D_1<\D_0$ and finite sets 
$s_0,s_1,\ldots,s_n$ such that $s_0<s_1<\cdots < s_n$ and $(s,\D_i)\in\G$ 
for every $s\in FU((s_1,\ldots,s_{i-1}))$ and $1\le i\le n$.

We will construct $\D_{n+1}$ and $s_{n+1}$. 
Let $\{ t_1,\ldots, t_k\} = FU ((s_0,s_1,\ldots, s_n))$. 
According to condition (i), there exist infinite disjoint collections 
$\D_{n+1}^1,\ldots, \D_{n+1}^k$ such that \newline 
$\D_{n+1}^k <\cdots < \D_{n+1}^1 < \D_n -s_n$ and 
$(t_i,\D_{n+1}^i)\in \G$ for every $1\le i\le k$. 
Set $\D_{n+1} = \D_{n+1}^k$ and choose $s_{n+1} \in \D_{n+1}$. 
Of course $s_n <s_{n+1}$. 
According to condition (ii), $(t_i,\D_{n+1})\in \G$ for all $1\le i\le k$. 

Set $\D = \{s_0,s_1,s_2,\ldots\}$ which is an infinite disjoint collection 
and $\D<\D_0$. 
Let $s\in FU (\D)$ and $\D' <\D-s$. 
If $s=\emptyset$, then $(\emptyset,\D_1)\in\G$ and according to (ii), 
$(\emptyset,\D')\in\G$.
If $s\ne \emptyset$, then $s\in FU((s_0,s_1,\ldots,s_{n_0}))$ for 
$n_0 = \min \{n\in\nat \cup \{0\} : s\in FU((s_0,s_1,\ldots,s_n))\}$. 
Thus $(s,\D_{n_0+1})\in \G$ and according to (ii), $(s,\D') \in\G$.
\end{proof}

We are now ready to state and prove the block Ramsey partition 
theorem for every countable ordinal number.

\begin{thm}
[Block-Ramsey partition theorem for every countable ordinal]
\label{thm:block-Ramsey}
Let \newline 
$\F\subseteq\B^{<\omega}(\nat)$ be a family of finite disjoint collections, 
$\D_0$ an infinite disjoint collection and $\xi$ a countable ordinal. 
Then there exists an infinite disjoint collection $\D<\D_0$ such that 
\begin{itemize}
\item[{}] either $\B^\xi \cap \B^{<\omega} (\D) \subseteq \F$ or 
$\B^\xi \cap \B^{<\omega}(\D)\subseteq\B^{<\omega} (\nat)\setminus \F$.
\end{itemize}
\end{thm}

\begin{proof} 
For $\xi=1$ it follows from Theorem~\ref{thm:hindman}.

Let $\xi>1$. 
Assume that the theorem is valid for every $\zeta <\xi$.
Let $s$ be a non-empty, finite subset of $\nat$ with $\min s=n$ and 
$\D$ an infinite disjoint collection with $\D<\D_0$. 
According to Proposition~\ref{justification}, there exist $\xi_n<\xi$ 
such that $\B^\xi (s) = \B^{\xi_n} \cap (\B^{<\omega} (\nat)-s)$. 
Using the induction hypothesis, there exists an infinite disjoint 
collection $\D_1 <\D$ such that 
\begin{itemize}
\item[{}] either $\B^{\xi_n} \cap \B^{<\omega}(\D_1) \subseteq \F(s)$;
\quad or\quad  
$\B^{\xi_n} \cap \B^{<\omega} (\D_1) \subseteq \B^{<\omega} (\nat)
\setminus \F(s)$. 
\end{itemize}
Set $\D_s = \D_1 - s<\D$. 
Then $s<\bigcup_{F\in \D_s} F$ and 
\begin{itemize}
\item[{}] either $\B^\xi (s) \cap \B^{<\omega} (\D_s) \subseteq \F(s)$;
\quad or \quad 
$\B^\xi (s) \cap \B^{<\omega} (\D_s) \subseteq \B^{<\omega} (\nat)
\setminus \F(s)$.
\end{itemize}

Let $\G= \{(s,\D): s\in [\nat]^{<\omega}$, $\D\in \B^\omega (\nat)$ 
with $\D<\D_0$, $s<\bigcup_{F\in\D}F$ and 
\begin{itemize}
\item[{}] either $\B^\xi (s)\cap \B^{<\omega} (\D) \subseteq \F(s)$;
\quad or\quad  
$\B^\xi (s)\cap \B^{<\omega} (\D) \subseteq \B^{<\omega} (\nat) 
\setminus \F(s)\}$.
\end{itemize}
The family $\G$ satisfies the conditions (i) (by the above arguments) and 
(ii) (obviously) of Lemma~\ref{lem:block-Ramsey}. 
hence there exists an infinite disjoint collection $\D_1<\D_0$ such that 
$(s,\D'_1)\in\G$ for every $s\in FU(\D_1)$ and $\D'_1 <\D_1 -s$.

Let $\F_1 = \{(s)\in FU(\D_1) :\B^\xi (s)\cap \B^{<\omega}
(\D_1-s) \subseteq \F(s)\}$. 
We use the induction hypothesis for $\xi=1$ (Theorem~\ref{thm:hindman}). 
Then there exists an infinite disjoint collection $\D<\D_1$ such that 
\begin{itemize}
\item[{}] either $\B^1\cap\B^{<\omega} (\D)\subseteq \F_1$;
\quad or \quad 
$\B^1\cap\B^{<\omega}(\D)\subseteq \B^{<\omega}(\nat)\setminus \F_1$. 
\end{itemize}
Since $\D<\D_1$ we have $\B^{<\omega} (\D)\subseteq \B^{<\omega}(\D_1)$ 
and consequently that $(s,\D-s) \in\G$ for every $s\in FU(\D)$. 
Thus 
\begin{itemize}
\item[{}] either $\B^\xi (s)\cap\B^{<\omega}(\D-s)\subseteq \F(s)$ for every 
$s\in FU(\D)$;  
\item[{}] or $\B^\xi (s)\cap\B^{<\omega}(\D-s)\subseteq \B^{<\omega}(\nat)
\setminus \F(s)$ for every $s\in FU(\D)$.
\end{itemize}
Hence, 
\begin{itemize}
\item[{}] either $\B^\xi\cap\B^{<\omega}(\D)\subseteq \F$;
\quad or\quad  
$\B^\xi\cap \B^{<\omega}(\D) \subseteq \B^{<\omega}(\nat)\setminus \F$.
\qquad \qed
\end{itemize}
\renewcommand{\qed}{}
\end{proof}

For finite $\xi (<\omega)$ our block Ramsey theorem is similar to 
Theorem~2.2 in Milliken \cite{M}, (stated there in terms of the 
corresponding sum-sets). 
%and also follows from Carlson's general partition theorem (\cite{C}). 
For infinite countable ordinals $\xi$ it is new.

For completeness' sake we now state the corresponding result for sum-sets.

\begin{defn}\label{def:sum-sets}
1. The natural correspondence between $\nat$ and $[\nat]^{<\omega}$ is 
given by the function 
$$\varphi :[\nat]^{<\omega} \to \nat\text{ with } \varphi (s) = 
\sum_{\ell\in s} 2^{\ell-1}\text{ for } s\in [\nat]^{<\omega} \ .$$
The function $\varphi$ has the following properties:
\begin{itemize}
\item[(i)] $\varphi$ is bijection;
\item[(ii)] if $s_1<s_2$, then $\varphi (s_1) <\varphi (s_2)$; and 
\item[(iii)] $\varphi (s_1)+\varphi (s_2)= \varphi (s_1\cup s_2)$ if 
$s_1\cap s_2=\emptyset$.
\end{itemize}

2. For every $n\in\nat$ we define $c(n) = \min \varphi^{-1}(n)$.

3. For every subset $L$ of $\nat$ the 1-sum set $\P^1(L)$ of $L$ is 
defined as follows:
$$\P^1 (L) = \bigg\{\sum_{\ell\in s} \ell : s\in [L]^{<\omega}\bigg\}\ .$$

4. For every countable ordinal $1\le \xi<\omega_1$, we define the $\xi$-sum 
set $\P^\xi(L)$ for every infinite subset $L$ of $\nat$ as follows: 
\begin{equation*}
\begin{split}
\P^\xi (L) =
& \bigg\{\bigg(\sum_{\ell\in s_1}\ell,\ldots,\sum_{\ell\in s_k}\ell\bigg) :
k\in\nat,\ s_1<\cdots < s_k \in [L]^{<\omega} \\
&\qquad \text{ and }\ (c(\min s_1 ),\ldots, c(\min s_k)) \in \A_\xi\bigg\}\ .
\end{split}
\end{equation*}

5.  For every $k\in\nat$ and $L\in [\nat]$
\begin{equation*}
\begin{split}
&\P^k(L) =\bigg\{\bigg(\sum_{\ell\in s_1} \ell,\ldots, 
\sum_{\ell\in s_k}\ell\bigg) : s_1 < \cdots < s_k\in [L]^{<\omega}\bigg\}\ 
\text{ and }\\ 
&\P(L) = \bigcup_{k\in\nat} \P^k (L)\ . 
\end{split}
\end{equation*}
\end{defn}

\begin{thm}[Hindman's sum-theorem extended to countable ordinals] 
\label{thm:hindmanextended}
Let $\F \subseteq [\nat]^{<\omega}$. 
Then there exists $L\in [\nat]$ such that either 
$\P^\xi (L) \subseteq \F$; or $\P^\xi (L)\subseteq [\nat]^{<\omega}
\setminus \F$.
\end{thm}

\begin{proof} 
Let $\bar\varphi :\B^{<\omega} (\nat) \to [\nat]^{<\omega}$ with 
$\bar \varphi ((s_1,\ldots, s_k)) = (\varphi (s_1),\ldots,\varphi(s_k))$, 
and let \break 
$\F_1 = \bar\varphi^{-1} (\F)\subseteq \B^{<\omega}(\nat)$. 

According to the $\xi$-block-Ramsey partition theorem 
(Theorem~\ref{thm:block-Ramsey}), 
there exists an 
infinite disjoint collection $\D$ such that 
\begin{itemize}
\item[{}] either $\B^\xi\cap \B^{<\omega}(\D)\subseteq \F_1$;
\quad or\quad  
$\B^\xi\cap \B^{<\omega} (\D) \subseteq \B^{<\omega} (\nat)\setminus \F_1$.
\end{itemize}
Equivalently, 
\begin{itemize}
\item[{}] either $\bar\varphi (\B^\xi\cap \B^{<\omega}(\D)\subset \F$; 
\quad or\quad  
$\bar \varphi (\B^\xi\cap \B^{<\omega}(\D)\subseteq [\nat]^{<\omega}
\setminus \F$.
\end{itemize}

Set $L = \{\varphi (t) :t\in\D\}$. 
We will prove that $\P^\xi (L)\subseteq \bar\varphi (\B^\xi \cap 
\B^{<\omega} (\D)$. 
Let 
$$\bigg(\sum_{\ell\in s_1}\ell,\ldots, \sum_{\ell\in s_k}\ell\bigg) 
\in \P^\xi(L)\ .$$ 
Then $s_1<\cdots < s_k\in [L]^{<\omega}$.
Thus for every $i=1,\ldots,k$, we have that $s_i = (\varphi (t_1^i),\ldots,
\varphi (t_{k_i}^i))$ for $t_1^i <\cdots < t_{k_i}^i\in\D$. 
Hence, $\sum_{\ell\in s_i}\ell = \varphi (t_1^i)+\cdots + \varphi (t_{k_i}^i) 
= \varphi (t_1^i\cup\cdots\cup t_{k_i}^i)$ for every $i\in I$. 
Set $t_i = t_1^i\cup\cdots\cup t_{k_i}^i\in FU(\D)$. 
Then $\sum_{\ell\in s_i} \ell = \varphi (t_i)$ and 
$t_i \in FU(\D)$ for  $1\le i\le k$. 
Hence $(\sum_{\ell\in s_1} \ell,\ldots,\sum_{\ell\in s_k}\ell) 
= \bar\varphi ((t_1,\ldots,t_k))$.
Since $\D$ is an infinite disjoint collection, for every $t,t'\in\D$ we have 
that $t<t'$ if $\varphi (t) <\varphi (t')$. 
Hence, since $s_1<\cdots < s_k$, we have that $t_1<\cdots <t_k$. 
Thus $(t_1,\ldots, t_k)\in\B^{<\omega} (\D)$. 
Since $(c(\min s_1),\ldots, c(\min s_k))\in \A_\xi$  and 
$c(\min s_i) = \min t_i$ for all $1\le i\le k$
we have that $(t_1,\ldots, t_k)\in \B^\xi\cap \B^{<\omega}(\D)$ and 
consequently that 
$$\bigg(\sum_{\ell\in s_1} \ell,\ldots, \sum_{\ell\in s_k}\ell\bigg) 
\in \bar\varphi (\B^\xi \cap \B^{<\omega} (\D))\ .$$
Thus $\P^\xi(L) \subseteq \bar\varphi (\B^\xi\cap \B^{<\omega}(\D))$. 
Hence, 
either $\P^\xi (L) \subseteq \F$; or 
$\P^\xi (L) \subseteq [\nat]^{<\omega} \setminus \F$.
\end{proof}

\begin{cor}\label{cor:colored1}
If $[\nat]^{<\omega}$ is finitely colored, then for every countable 
ordinal $\xi$ there exist an infinite subset $L$ of $\nat$ such that 
the set $\P^\xi(L)$ is monochromatic.
\end{cor}

Especially, for the case $\xi=\omega$ the statement of 
Corollary~\ref{cor:colored1} is the following: 

\begin{cor}\label{cor:colored2}
If $[\nat]^{<\omega}$ is finitely colored, then there exists $L\in[\nat]$ 
such that the set $\P^\omega (L) = \{(\sum_{\ell\in s_1}\ell,\ldots, 
\sum_{\ell\in s_k}\ell) : k\in \nat$, $s_1<\cdots<s_k\in [L]^{<\omega}$ 
and $c(\min s_1)=k\}$ is monochromatic.
\end{cor} 

%%%%%%% 22222222
\section{Some properties of the recursive block Schreier system}

The main results proved in this section concerns 
the thinness of the families $\B^\xi$ (Proposition~\ref{prop:thinfamily}), 
the canonical representation of every non-empty disjoint collection with 
respect to $\B^\xi$ (Proposition~\ref{prop:canonicalrep}), 
and the computation of the strong Cantor-Bendixson index of $\B^\xi$ 
(Proposition~\ref{prop:Cantor-Bendix}). 
These properties, and their proofs, for $\B^\xi$, are analogous to those 
for $\A_\xi$ in \cite{F3}; 
they are necessary for the proof of the main results in Section~3.

\begin{defn}\label{def:Fthin}
Let $\F\subseteq \B^{<\omega}(\nat)$ be a family of finite disjoint 
collections.
\begin{itemize}
\item[(i)] $\F$ is {\em thin\/} if there are no elements $\bs,\bt\in\F$
with $\bs\ne \bt$ and $\bs \propto \bt$.
\item[(ii)] $\F^* = \{\bt \in \B^{<\omega}(\nat): \bt\propto \bs$ for some 
$\bs\in \F\}\cup \{\emptyset\}$.
\item[(iii)] $\F_* = \{\bt \in \B^{<\omega} (\nat): \bt\subseteq FU(\bs)$ 
for some $\bs\in \F\}\cup\{\emptyset\}$.
\item[(iv)] $\F$ is a {\em tree\/} if $\F^* = \F$.
\item[(v)] $\F$ is {\em hereditary\/} if $\F_* = \F$.
\end{itemize}
\end{defn}

Thinness and canonical representation are two properties of the Schreier
system $(\A_\xi)_{\xi<\omega_1}$ that were proved useful in \cite{F3} 
(see Proposition~\ref{prop:recursive}(ii) and (iii)). 
We will now prove that the block Schreier system $(\B^\xi)_{\xi<\omega_1}$ 
satisfies these properties too.

\begin{prop}\label{prop:thinfamily}
Every family $\B^\xi$, for $\xi<\omega_1$ is thin.
\end{prop} 

\begin{proof} 
The family $\B^1 = \{(s) :s\in [\nat]^{<\omega}\}$ is obviously thin.
Let $\xi>1$. 
Then ${\emptyset \notin \B^\xi}$. 
Let $\bs = (s_1,\ldots,s_k)\in \B^\xi$, $\bt =(t_1,\ldots,t_\lambda)\in\B^\xi$
such that $\bs \propto \bt$.
Then ${s = (\min s_1,\ldots,\min s_k)\in \A_\xi}$, 
$t= (\min t_1,\ldots,\min t_\lambda)\in \A_\xi$ and $s\propto t$. 
Since $\A_\xi$ is thin (Proposition~\ref{prop:recursive}(ii)), we have 
that $s=t$ and consequently that $\bs = \bt$. 
Hence, $\B^\xi$ is thin.
\end{proof}

In the following we will prove that every disjoint collection (finite or 
infinite) has unique canonical representation with respect to each family 
$\B^\xi$.

\begin{defn}\label{def:non-zero}
Let $\xi$ be a non-zero countable ordinal number.

(i) A non-empty finite disjoint collection $\bs$ has canonical representation 
with respect to $\B^\xi$ if there exist unique $n\in\nat$, 
$\bs_1,\ldots,\bs_n\in \B^\xi$ and a, possibly empty, set $\bs_{n+1}\in 
(\B^\xi)^* \setminus \B^\xi$ with 
$\bigcup_{F\in \bs_1} F <\cdots < \bigcup_{F\in \bs_{n+1}}F$ and 
$\bs = \bigcup_{i=1}^{n+1}\bs_i$.

(ii) An infinite disjoint collection $\D$ has canonical representation with 
respect to $\B^\xi$ if there exists unique sequence $(\bs_n)_{n\in\nat}
\subseteq \B^\xi$ with $\bigcup_{F\in\bs_1} F<\bigcup_{F\in\bs_2} F<\cdots$ 
and ${\D= \bigcup_{n\in\nat} \bs_n}$. 
\end{defn}

\begin{prop}\label{prop:canonicalrep}
Let $\xi$ be a non-zero countable ordinal number. 
Every non-empty disjoint collection (finite or infinite) has canonical 
representation with respect to $\B^\xi$.
\end{prop} 

\begin{proof}
(i) Let $\D = (\bt_n)_{n\in\nat}$ be an infinite disjoint collection and 
$m_n=\min t_n$ for $n\in\nat$. 
Of course, $m_n<m_{n+1}$. 
The infinite subset $I= \{m_n :n\in\nat\}$ of $\nat$ has canonical 
representation with respect to $\A_\xi$ 
(Proposition~\ref{prop:recursive}(iii)).
Hence, there exists a unique sequence $(s_n)_{n\in\nat}\subseteq \A_\xi$ 
such that $I = \bigcup_{n=1}^\infty s_n$ and $s_1<s_2<\cdots$.
For every $n\in\nat$, set $\bs_n = \{t_i: i\in \nat$ and $m_i\in s_n\}$. 
Then $(\bs_n)$ is a canonical representation of $\D$ with respect to 
$\B^\xi$, and it is unique, since $\B^\xi$ is thin.

(ii) Let $\bs = (t_1,\ldots,t_k)$ be a finite disjoint collection. 
Set $t_{k+i} = \{n_k +i\}$ for every $i\in\nat$, where $n_k=\max t_k$. 
The infinite disjoint collection $\D = (t_n)_{n\in\nat}$ has canonical 
representation with respect to $\B^\xi$, according to (i). 
Using this fact, it is easy to prove that $\bs$ has canonical 
representation with respect to $\B^\xi$. 
It is unique, since $\B^\xi$ is thin.
\end{proof}

Identifying every finite collection $\bs$ 
(resp.\ every infinite collection $\D$) 
of finite subsets of $\nat$ with its 
characteristic function $x_{\bs} \in \{0,1\}^{[\nat]^{<\omega}}$ 
(resp.\ $x_{\D}\in \{0,1\}^{[\nat]^{<\omega}}$), 
we topologize the set of all the collections of finite subsets of 
$\nat$ by the topology of pointwise convergence (equivalently 
by the product topology of $\{0,1\}^{[\nat]^{<\omega}}$). 
Hence, we will say that a family $\F\in \B^{<\omega} (\nat)$ is pointwise
closed iff the family $\{x_{\bs} :\bs \in \F\}$ is closed in the 
topology of pointwise convergence in $\{0,1\}^{[\nat]^{<\omega}}$ and 
also a family $\U\subseteq \B^\omega (\nat)$ is pointwise closed iff 
$\{x_{\D} :\D\in \U\}$ is pointwise closed in $\{0,1\}^{[\nat]^{<\omega}}$.

We next turn our attention to trees and hereditary families. 
Closure of a tree or a hereditary family (in the relative topology 
of $\B^{<\omega}(\nat)$) is described by the non-existence of infinite
chains or infinite disjoint collections, respectively. 
The hereditary family generated by $\B^\xi$ (in any infinite disjoint 
collection) turns out to be relatively pointwise closed. 

\begin{prop}\label{prop:finitefamily}
Let $\F\subseteq \B^{<\omega}(\nat)$ be a family of finite disjoint 
collections.

{\rm (i)} If $\F$ is a tree, then $\F$ is pointwise closed if and only 
if there does not contain an infinite sequence $(\bs_n)_{n\in\nat}$ 
such that $\bs_1 \propto \bs_2\propto \cdots$ and $\bs_n\ne \bs_{n+1}$ 
for all $n\in\nat$.

{\rm (ii)} If $\F$ is hereditary, then $\F$ is pointwise closed if 
and only if there does not exist an infinite disjoint collection 
$\D$ such that $\B^{<\omega} (\D) \subseteq \F$.
\end{prop}

\begin{proof}
(i) Let $\F$ be a tree. 
If $(\bs_n)_{n\in\nat} \subseteq \F$, $\bs_1 \propto \bs_2\propto\cdots$ 
and $\bs_n\ne \bs_{n+1}$ for all $n\in\nat$, then $(\bs_n)$ converges 
in the product topology of $\{0,1\}^{[\nat]^{<\omega}}$ to the infinite 
disjoint collection $\D = \bigcup_{n=1}^\infty \bs_n$, which does not 
belong to $\F$, hence $\F$ is not pointwise closed. 

We assume that there does not exist a sequence $(\bs_n)\subseteq \F$ 
with $\bs_1 \propto \bs_2 \propto\cdots $ and $\bs_n\ne \bs_{n+1}$ 
for all $n\in\nat$. 
We will prove that $\F$ is pointwise closed. 
Let $(\bt_i)_{i\in\nat}\subseteq \F$ which converges pointwise 
to some $\bt\subseteq [\nat]^{<\omega}$. 
Let $t_1,t_2\in \bt$. 
Then, there exists $\iota_0\in\nat$ such that $\{t_1,t_2\}\subseteq 
t_{\iota_0}$, hence either $t_1<t_2$ or $t_2<t_1$. 
Thus $\bt$ is a disjoint collection.
Let $\bt$ be finite.
Then $\bt = (t_1<\cdots < t_k)$ for some $k\in\nat$ or $\bt =\emptyset$.
If $\bt=\emptyset$, then $\bt \in \F$ since $\F$ is a tree. 
Let $\bt = (t_1,\ldots, t_k)$. 
Since $FU ((\{1\},\ldots, \{\max t_k\}))$ is a finite set, there exists 
$\iota_0\in\nat$ such that $\bt \propto t_{\iota_0}$. 
Since $\F$ is a tree, $\bt \in\F$. 

We will prove that $\bt$ can not be an infinite set. 
Let $\bt = (t_n)_{n\in\nat}$ be an infinite disjoint collection. 
Set $\bs_n = (t_1 <\cdots < t_n)$ for every $n\in\nat$. 
Then $\bs_1 \propto \bs_2\propto \cdots $ and $\bs_n \ne \bs_{n+1}$ 
for all $n\in\nat$. 
We will prove that $(\bs_n)\subseteq\F$, which is a contradiction 
to our assumption. 
Indeed, let $n\in\nat$ and $\bt_i^n = \bt_i\cap FU((\{1\},\ldots,
\{\max t_n\}))$ for $i\in \nat$. 
Since $\bt_i^n \propto \bt_i$ and $\F$ is a tree we have that 
$\bt_i^n\in\F$ for all $n,i\in\nat$. 
The sequence $(\bt_i^n)_{i\in\nat}\subseteq \F$ converges pointwise to 
$\bs_n = (t_1<\cdots < t_n)$, since $\bt_i\to \bt$ and $FU((\{1\},\ldots,
\{\max t_n\}))$ is finite. 
Since $\bs_n$ is finite we have that $\bs_n\in\F$. 
Hence $\F$ is pointwise closed.

(ii) Let $\F$ be hereditary. 
If there exists an infinite disjoint collection $\D = (s_n)_{n\in\nat}
\subseteq \F$ such that $\B^{<\omega} (\D)\subseteq \F$, then 
$\bs_n = (s_1<\cdots < s_n)\in \F$ for every $n\in \nat$. 
According to case (i), $\F$ is not pointwise closed, since it is 
a tree. 
On the other hand,
if $\F$ is not pointwise closed, according to (i), there exist $\D= 
(s_n)_{n\in\nat}$ such that $(s_1,\ldots ,s_n)\in\F$ for every $n\in\nat$. 
Then $\B^{<\omega} (\D)\subseteq \F$. 
Indeed, if $\bt \in \B^{<\omega} (\D)$, 
then there exists $n\in\nat$ such that\newline 
$\bt\subseteq FU((s_1,\ldots, s_n)) \subseteq \F$, since $\F$ is hereditary.
\end{proof}

\begin{cor}\label{cor:pointwise}
Let $\F\subseteq \B^{<\omega} (\nat)$ which  is hereditary and pointwise 
closed. 
Then, every hereditary family $\F_1$ with 
$\F_1 \subseteq \F$ is also pointwise closed.
\end{cor}

\begin{cor}\label{cor:countableordinal}
Let $\xi$ be a countable ordinal and $\D$ an infinite disjoint collection 
Then the hereditary family $(\B^\xi \cap \B^{<\omega} (\D))_*$ 
is pointwise closed.
\end{cor} 

\begin{proof} 
The family $(\B^1 \cap \B^{<\omega} (\D))_*$ is obviously pointwise closed. 
Let $\xi>1$, and assume that $(\B^\zeta \cap \B^{<\omega} (\D))_*$ 
is pointwise closed for every $\zeta <\xi$. 
Let $(\B^\xi \cap \B^{<\omega} (\D))_*$ is not pointwise closed. 
According to Proposition~\ref{prop:finitefamily}, there exists an 
infinite disjoint collection $\D_0 = (s_n)_{n\in\nat}$ such that 
$\B^{<\omega} (\D_0) \subseteq (\B^\xi \cap \B^{<\omega} (\D))_*$. 
%Let $k_0 = \min s_1$ and $n\in\nat$. 
Let $n\in \nat$.
Then, $\bs_n = (s_1 ,\ldots, s_n) \in (\B^\xi \cap \B^{<\omega}(\D))_*$,  
so there exists $\bt = (t_1,\ldots, t_\lambda)\in \B^\xi \cap 
\B^{<\omega}(\D)$ such that $\bs_n \subseteq FU(\bt)$. 
If $k= \min t_1$, then  according to 
Proposition~\ref{justification}, $(t_2,\ldots,t_\lambda)\in 
\B^{\xi_k} \cap \B^{<\omega} (\D-t_1)$. 
Hence, 
$(s_2,\ldots, s_n) \in (\B^{\xi_k} \cap \B^{<\omega} (\D))_*$ with 
$k\le \min s_1$. 
Thus, there exists $k_1\le \min s_1$ such that $(s_2,\ldots,s_n) \in 
(\B^{\xi_{k_1}} \cap \B^{<\omega} (\D))_*$ for infinitely many $n\in\nat$. 
This gives that $(s_2,\ldots,s_n)\in (\B^{\xi_{k_1}} \cap \B^{<\omega} 
(\D))_*$ for every $n\in\nat$, and consequently 
that $(\B^{\xi_{k_1}}\cap \B^{<\omega}(\D))_*$ 
is not pointwise closed,  according to Proposition~\ref{prop:finitefamily}.
This is a contradiction to the induction hypothesis, since $\xi_{k_1} <\xi$.
\end{proof}

For hereditary and pointwise closed families $\F$ in $\B^{<\omega}(\nat)$ 
we can define the strong Cantor-Bendixson index $s_\D(\F)$ with 
respect to an infinite disjoint collection $\D$. 
The strong Cantor-Bendixson index for $\B^\xi$ will turn out to be 
$\xi+1$, for any infinite disjoint collection.

\begin{defn}\label{def:Cantor-Bendix}
Let $\F\subseteq \B^{<\omega}(\nat)$ be a hereditary and pointwise closed 
family of finite disjoint collections. 
For an infinite disjoint collection $\D$ we define the {\em strong block 
Cantor-Bendixson derivatives} $(\F)_{\D}^\xi$ of $\F$ on $\D$ for every 
$\xi <\omega_1$ as follows: 
\begin{equation*}
\begin{split}
(\F)_\D^\xi = &\{\bs \in \F\cap \B^{<\omega} (\D) :\text{ the set }
\{t\in FU(\D) :\bs \cup \{t\}\notin \F\} \\
\noalign{\vskip6pt}
&\text{does not contain an infinite disjoint collection}\}.
\end{split}
\end{equation*}
It is easy to verify that $(\F)_\D^1$ is hereditary, hence pointwise closed.
So, we can define for every $\xi >1$ the $\xi$-derivatives of $\F$ 
recursively as follows:
\begin{gather*} 
(\F)_\D^{\zeta +1} = ((\F)_\D^\zeta  )_\D^1\ \text{ for all } \zeta <\omega_1;
\ \text{ and}\\
(\F)_\D^\xi = \bigcap_{\beta<\xi} (\F)_\D^\beta\ 
\text{ for $\xi$ a limit ordinal}\ .
\end{gather*}

The  {\em strong block Cantor-Bendixson index} $s_\D(\F)$ of $\F$ on $\D$ 
is the smallest countable ordinal $\xi$ such that $(\F)_\D^\xi = \emptyset$.
\end{defn}

\begin{remark}\label{rem:Cantor-Bendix}
{\rm (i)} The strong block Cantor-Bendixson index $s_\D(\F)$ of a hereditary 
and pointwise closed family of finite disjoint collections on an infinite
disjoint collection $\D$ is a countable successor ordinal less than or 
equal to the ``usual'' Cantor-Bendixson index $O(\F)$ of $\F$ into 
$\{0,1\}^{[\nat]^{<\omega}}$ (see \cite{K}).
\begin{itemize}
\item[(ii)] $s_\D (\F\cap \B^{<\omega} (\D)) = s_\D (\F)$.

\item[(iii)] $s_\D(\F_1) \le s_\D (\F_2)$ if $\F_1,\F_2 \subseteq \B^{<\omega} 
(\nat)$ are hereditary and pointwise closed families with $\F_1\subseteq 
\F_2$.

\item[(iv)] If $\bs \in (\F)_\D^\xi$, then $\bs \cap FU (\D_1)\in 
(\F)_{\D_1}^\xi$ for every $\D_1 <\D$, since $FU (\D_1)\subseteq FU (\D)$.

\item[(v)] If $\D_1<\D$, then $s_{\D_1} (\F) \ge s_\D  (\F)$, 
according to (iv).

\item[(vi)] If $\D_1\setminus \D$ is a finite set, then $s_{\D_1}(\F) 
\ge s_\D(\F)$.
\end{itemize}
\end{remark}

\begin{prop} %[strong block Cantor-Bendixson index for $\B^\xi$] 
\label{prop:Cantor-Bendix}
Let $\D$ be an infinite disjoint collection and $\D_1<\D$. 
Then 
$$s_{\D_1}\Big((\B^\xi \cap \B^{<\omega} (\D))_*\Big) 
= \xi+1\ \text{ for every }\ \xi<\omega_1\ .$$
\end{prop}

\begin{proof} 
For every $\xi<\omega_1$ the families $(\B^\xi \cap \B^{<\omega} (\D))_*$ 
are pointwise closed (Corollary~\ref{cor:countableordinal}). 
Also, for every $s\in FU(\D)$ with $\min s=n$ we have, according to 
Proposition~\ref{justification}, that 
$$(\B^\xi \cap \B^{<\omega}(\D)) (s) = \B^{\xi_n} \cap \B^{<\omega} 
(\D-s)\text{ for some } \xi_n <\xi\ .$$

We will prove by induction that 
$((\B^\xi\cap \B^{<\omega}(\D))_*)_{\D_1}^\xi 
= \{\emptyset\}\ \text{ for every }\ \xi<\omega_1$,
Of course, 
$(\B^1\cap \B^{<\omega} (\D))_* = \{\{s\} :s\in FU(\D)\} \cup \{\emptyset\}$.
hence $((\B^1 \cap \B^{<\omega} (\D))_*)_{\D_1}^1 = \{\emptyset\}$.
%and consequently 
%$$s_{\D_1} ((\B^1 \cap \B^{<\omega} (\D))_*) = 2\ .$$

Let $\xi>1$ and assume that 
$$((\B^\zeta \cap \B^{<\omega} (\D))_* )_{\D_1}^\zeta 
= \{\emptyset\}\text{ for every }\zeta <\xi\text{ and } \D_1<\D\ .$$
Hence, for every $s\in FU(\D_1)$ with $\min s=n$ we have that 
$$((\B^\xi \cap \B^{<\omega} (\D))(s)_*)_{\D_1}^{\xi_n} 
= ((\B^{\xi_n} \cap  \B^{<\omega} (\D-s))_*)_{\D_1}^{\xi_n} 
= \{\emptyset\} \ .$$
This gives that 
$\{(s)\} \in ((\B^\xi\cap \B^{<\omega} (\D))_*))_{\D_1}^{\xi_n}$.
So, if $\xi = \zeta+1$ is a successor ordinal, we have that 
$\{(s)\} \in ((\B^\xi\cap \B^{<\omega}(\D))_*)_{\D_1}^\xi$  
for every $s\in FU(\D_1)$, 
hence $\emptyset \in ((\B^\xi \cap \B^{<\omega}(\D))_*)_{\D_1}^\xi$ 
and if 
$\xi$ is a limit ordinal, since $\emptyset\in ((\B^\xi\cap \B^{<\omega} 
(\D))_*)_{\D_1}^{\xi_n}$ for every $n\in\nat$ 
and $\sup \xi_n  =\xi$, also we have that 
$\emptyset \in ((\B^\xi \cap \B^{<\omega} (\D))_*)_{\D_1}^\xi$. 

If $\{\emptyset\} \ne ((\B^\xi \cap \B^{<\omega} (\D))_*)_{\D_1}^\xi$, 
then there exist $D_2<D_2$ and  $s\in FU(\D_2)$ such that 
${((\B^\xi \cap \B^{<\omega} (\D))(s)_*)_{\D_2}^\xi \ne\emptyset}$ 
(see Lemma 2.8 in \cite{F3}). 
A  contradiction to the induction hypothesis. 
Hence, $((\B^\xi \cap \B^{<\omega}(\D))_*)_{\D_1}^\xi = \emptyset$.
\end{proof}

In case the partition family is a tree, we obtain, using the canonical 
representation, an alternative description for the second horn of the block
Ramsey dichotomy (Theorem~\ref{thm:block-Ramsey}). 

\begin{prop}\label{prop:tree}
Let $\F\subseteq \B^{<\omega}(\nat)$ be a tree, $\D$ an infinite disjoint 
collection and $\xi<\omega_1$. 
Then $\B^\xi\cap \B^{<\omega}(\D)\subseteq \B^{<\omega}(\nat)\setminus \F$
if and only if 
$\F\cap \B^{<\omega} (\D) \subseteq (\B^\xi)^* \setminus \B^\xi$.
\end{prop}

\begin{proof} 
Let $\B^\xi \cap \B^{<\omega} (\D) \subseteq \B^{<\omega} (\nat)\setminus \F$
and $\bs = (s_1,\ldots,s_k)\in \F\cap \B^{<\omega} (\D)$. 
Then $\bs$ has canonical representation with respect to $\B^\xi$
(Proposition~\ref{prop:canonicalrep}), hence 
either $\bs \in (\B^\xi)^*\setminus \B^\xi$, as required, or there 
exists $\bs_1\in \B_\xi$ such that $\bs_1 \propto \bs$. 
The second case is impossible. 
Indeed, since $\F$ is a tree and $\bs\in \F\cap \B^{<\omega} (\D)$, we have 
$\bs_1 \in \F\cap \B^{<\omega} (\D) \cap \B_\xi$; a contradiction to our 
assumption. 
Hence, $\F\cap \B^{<\omega} (\D) \subseteq (\B^\xi)^*\setminus \B^\xi$.

It is obvious that 
$\B^\xi \cap \B^{<\omega} (\D) \subseteq \B^{<\omega} (\nat)\setminus\F$
if $\F\cap \B^{<\omega}(\D) \subseteq (\B^\xi)^* 
\setminus \B^\xi$. 
\end{proof}

\begin{cor}\label{cor:tree} 
Let $\xi_1,\xi_2$ be countable ordinals with $\xi_1<\xi_2$. 
For every infinite disjoint collection $\D$ there exists an 
infinite disjoint collection $\D_1<\D$ such that 
\begin{equation*}
(\B^{\xi_1})_* \cap \B^{<\omega} (\D_1) \subseteq (\B^{\xi_2})^* 
\setminus \B^{\xi_2}\ .
\end{equation*}
\end{cor}

\begin{proof} 
Of course $(\B^{\xi_1})_*$ is a tree. 
According to the $\xi_2$-block Ramsey partition theorem 
(Theorem~\ref{thm:block-Ramsey}) and the previous 
proposition, for every infinite disjoint collection $\D$ there exists 
$\D_1 <\D$ such that either 
$\B^{\xi_2} \cap \B^{<\omega} (\D_1) \subseteq (\B^{\xi_1})_*$ or 
${(\B^{\xi_1})_* \cap \B^{<\omega} (\D_1) \subseteq (\B^{\xi_2})^*\setminus 
\B^{\xi_2}}$. 
The first alternative is impossible, since in this case, 
according to Proposition~\ref{prop:Cantor-Bendix}, 
$\xi_2 +1 = s_{\D_1} ((\B^{\xi_2} \cap \B^{<\omega} (\D_1))_*) 
\le s_{\D_1} ((\B^{\xi_1})_*) = \xi_1 +1$; 
a contradiction. 
\end{proof}

\section{Block Nash-Williams partition theorems for every 
countable ordinal}		%% section 3 

Let us recall that, according to the block Ramsey theorem 
(Theorem~\ref{thm:block-Ramsey}), for every countable ordinal $\xi$ and 
every partition of $\F$ of $\B^{<\omega}(\nat)$, 
there is an infinite disjoint collection $\D$, 
all of whose $\B^\xi$-finite blocks are either in the partition family $\F$ 
itself or in the complement $\F^c$; but Theorem~\ref{thm:block-Ramsey} 
provides no information on which to decide, whether in $\F$ or in $\F^c$. 
We now have at our disposal all the tools that will allow us to obtain 
for a partition family $\F$ that is hereditary 
(in Theorem~\ref{block-NashWilliams})  
and a tree (in Theorem~\ref{block-NashWilliams2}), a criterion, in 
terms of the strong Cantor-Bendixson index of $\F$, according to which, if 
this index is greater than $\xi +1$, all $\B^\xi$-finite blocks fall in $\F$, 
and if less than $\xi+1$, in $\F^c$ (albeit in a weaker, non-symmetrical 
manner). 

%We now have at our disposal the tools that will allow us to obtain 
%(in Theorem~\ref{block-NashWilliams}) for partitions that are hereditary 
%families, a very specific description of whether the block-Ramsey dichotomy 
%of Theorem~\ref{thm:block-Ramsey} will be resolved in favor of the 
%complement of $\F$. 
%This criterion, given in terms of the strong Cantor-Bendixson index of $\F$, 
%assigns to $\F$ the cases where the index of $\F$ is greater than $\xi+1$, 
%and $\F$ partitions the family $\B^\xi$, and to the complement of  $\F$, 
%in case it is less.

It will be observed that the dichotomy of Theorem~\ref{block-NashWilliams} 
is non-symmetric, reflecting the fact that the strong property of 
hereditariness is assumed for the family $\F$ itself; this type of 
non-symmetric dichotomies are characteristic of Nash-Williams 
partition theorems, and in fact, from Theorem~\ref{block-NashWilliams}, we 
will derive in the sequel 
(Theorem~\ref{block-NashWilliams2}, 
Corollaries~\ref{cor:tree2}, \ref{cor:sumsets} and \ref{cor:gowers}) 
various forms of block Nash-Williams theorems. 

\begin{thm}
%[Block Nash-Williams partition theorem for every countable ordinal for 
%hereditary families]
\label{block-NashWilliams}
Let $\F$ be a hereditary family of finite disjoint collections on $\nat$ 
and $\D$ an infinite disjoint collection on $\nat$. 
We have the following cases:

\noindent {\bf [Case 1]}
The family $\F\cap \B^{<\omega} (\D)$ is not pointwise closed. 
Then, there exists an infinite disjoint collection $\D_1<\D$ such that 
$\B^{<\omega} (\D_1) \subseteq \F$.

\noindent {\bf [Case 2]}
The family $\F\cap \B^{<\omega} (\D)$ is pointwise closed. 
Then, setting
\begin{equation*}
\xi_{\D}^{\F} = \sup \{s_{\D_1} (\F) : \D_1 <\D\}\ ,
\end{equation*}
which is a countable ordinal, the following subcases obtain:
\begin{itemize}
\item[2(i)] If $\xi+1 <\xi_{\D}^{\F}$, then there exists $\D_2 <\D$ 
such that 
\begin{equation*}
\B^\xi \cap \B^{<\omega} (\D_2) \subseteq \F\ ;
\end{equation*}
\item[2(ii)] if $\xi+1>\xi_{\D}^{\F}$, then for every $\D_1 <\D$ there exists 
$\D_2 <\D_1$ such that 
\begin{equation*}
\B^\xi \cap \B^{<\omega} (\D_2)\subseteq \B^{<\omega} (\nat)\setminus \F\ ;
\end{equation*}
(equivalently $\F\cap \B^{<\omega} (\D_2)\subseteq (\B^\xi)^*\setminus 
\B^\xi$); and  
\item[2(iii)] if $\xi+1 = \xi_{\D}^{\F}$, then there  exists $\D_2<\D$ 
such that 
\item[{}] either $\B^\xi \cap \B^{<\omega} (\D_2)\subseteq \F\ \text{ or }\  
\B^\xi \cap \B^{<\omega} (\D_2)\subseteq \B^{<\omega} (\nat)\setminus\F\ .$ 
\item[{}] Both alternatives may materialize.
\end{itemize}
\end{thm}

\begin{proof}{} 
[Case 1] If the hereditary family $\F\cap \B^{<\omega} (\D)$ is not pointwise 
closed, then, according to Proposition~\ref{prop:finitefamily}, there exists 
$\D_1<\D$ such that $\B^{<\omega}(\D_1) \subseteq \F\cap \B^{<\omega}(\D) 
\subseteq \F$.

\noindent [Case 2] 
If the hereditary family $\F\cap \B^{<\omega} (\D)$ is pointwise closed, 
then $\xi_{\D}^{\F}$ is countable, since the ``usual'' Cantor-Bendixson 
index $O(\F)$ of $\F$ into $\{0,1\}^{[\nat]^{<\omega}}$ is countable 
%(Remark~\ref{rem:block-NashWilliams}(i)) 
(Remark 2.9(i))
and $s_{\D_1}(\F)\le O(\F)$ for every 
$\D_1 <\D$.

2(i) 
Let $\xi+1 <\xi_{\D}^{\F}$. 
Then there exists $\D_1<\D$ such that $\xi+1 <s_{\D_1}(\F)$. 
According to Theorem~\ref{thm:block-Ramsey} and Proposition~\ref{prop:tree}, 
there exists $\D_2<\D_1$ such that 
\begin{itemize}
\item[{}] either $\B^\xi\cap \B^{<\omega}(\D_2)\subseteq \F$ 
\quad or\quad  
$\F\cap \B^{<\omega} (\D_2) \subseteq (\B^\xi)^*\setminus \B^\xi 
\subseteq (\B^\xi)_*$
\end{itemize}
The second alternative is impossible. 
Indeed, if $\F\cap \B^{<\omega}(\D_2)\subseteq (\B^\xi)_*$, then, according 
to Proposition~\ref{prop:Cantor-Bendix}, $s_{\D_2} (\F) = s_{\D_2}
(\F\cap \B^{<\omega} (\D_2)) \le s_{\D_2}((\B^\xi)_*) = \xi+1$; 
a contradiction. 
Hence, $\B^\xi \cap \B^{<\omega} (\D_2) \subseteq \F$.

2(ii) 
Let $\xi +1>\xi_{\D}^{\F}$, and $\D_1<\D$. 
According to the $\xi$-block Ramsey partition theorem
(Theorem~\ref{thm:block-Ramsey}), there exists 
$\D_2 <\D_1$ such that 
\begin{itemize}
\item[{}] either $\B^\xi \cap \B^{<\omega} (\D_2) \subseteq \F$
\quad or\quad  
$\B^\xi \cap \B^{<\omega}(\D_2)\subseteq \B^{<\omega} (\nat)\setminus \F$.
\end{itemize}
The first alternative is impossible, 
Indeed, if $(\B^\xi\cap \B^{<\omega} (\D_2)) \subseteq \F$, then, according 
to 
%Proposition~\ref{prop3.10} and Remark~\ref{rem:block-NashWilliams}, 
Proposition 2.10 and Remark 2.9, 
we obtain that 
\begin{equation*}
\xi+1 = s_{\D_2} ((\B^\xi \cap \B^{<\omega} (\D_2))_*) \le s_{\D_2} 
(\F) \le \xi_{\D}^{\F}\ ;
\end{equation*}
a contradiction. 
Hence, $\B^\xi\cap \B^{<\omega} (\D_2) \subseteq \B^{<\omega}(\nat) 
\setminus \F$ and according to Proposition~\ref{prop:tree}, 
$\F\cap \B^{<\omega} (\D_2)\subseteq (\B^\xi)^* \setminus \B^\xi$. 

2(iii) 
Let $\xi+1 = \xi_{\D}^{\F}$. 
According to Theorem~\ref{thm:block-Ramsey}, there exists $\D_2<\D$ 
such that 
\begin{itemize}
\item[{}] either $\B^\xi\cap \B^{<\omega} (\D_2)\subseteq \F$ or 
$\B^\xi\cap \B^{<\omega} (\D_2)\subseteq \B^{<\omega}(\nat)\setminus \F$.
\end{itemize}

That both alternatives may materialize can be seen by considering two 
simple examples: 

(1) $\F = \{\bs = (s_1<\cdots < s_{2k+1}) \in \B^{<\omega} (\nat): 
k\in\nat$ and $\min s_1 =k\}$. 
It is easy to see that $\F_*$ is pointwise closed (according to 
Proposition~\ref{prop:finitefamily}) and that	
${\F(s) = \B^{2n} \cap \B^{<\omega} (\nat-s)}$ for every $s\in [\nat]^{<\omega}$
with $\min s=n$. 
Analogously to Proposition~\ref{prop:Cantor-Bendix}, it can be proved that
$s_{\D}(\F_*) = \omega+1$ for every $\D\in \B^\omega (\nat)$.
It is now easy to verify that 
\begin{equation*}
\B^\omega \cap \B^{<\omega} (\D) \subseteq \F_*\ \text{ for every }\ 
\D\in \B^\omega (\nat)\ .
\end{equation*}

(2) $\F = \{\bs = (s_1 <\cdots < s_k) \in \B^{<\omega} (\nat):k\in\nat$ 
and $\min s_1 = 2k\}$. 
The family $\F_*$ is pointwise closed and $s_{\D_1}(\F_*) = \omega+1$ for 
every $\D_1 <\D$, where $\D$ is an infinite disjoint collection containing 
sets of even natural numbers. 
Thus, $\xi_{\D}^{\F} = \omega+1$. 
It is now easy to see that $\B^\xi\cap \B^{<\omega} (\D_1) \subseteq 
\B^{<\omega} (\nat) \setminus \F_*$ for every $\D_1 <\D$, since 
$\F\cap \B^{<\omega} (\D_1) \subseteq (\B^\xi)^* \setminus \B^\xi$.
\end{proof}

For the block Nash-Williams version of partition theorem for trees, given 
below (Theorem~\ref{block-NashWilliams2}) we need the following definition. 

\begin{defn}\label{def:blockNW} 
Let $\F$ be a family of finite disjoint collections on $\nat$. 
We set 
\begin{equation*}
\F_h = \{\bs \in \F:\text{ every }\bt \in FU(\bs)\text{ belongs to }\F\} 
\cup \{\emptyset\}.
\end{equation*}
Of course, $\F_h$ is the largest subfamily of $\F$ which is hereditary.
\end{defn}

\begin{thm}
%[Block Nash-Williams partition theorem for every countable ordinal for trees]
\label{block-NashWilliams2}
Let $\F$ be a family of finite disjoint collections on $\nat$ which is a tree 
and let $\D$ be an infinite disjoint collection on $\nat$. 
We have the following cases:

\noindent {\bf [Case 1]}
The family $\F_h\cap \B^{<\omega} (\D)$ is not pointwise closed.

Then, there exists an infinite disjoint collection $\D_1<\D$ such that 
$\B^{<\omega} (\D_1) \subseteq \F$.

\noindent {\bf [Case 2]}
The family $\F_h \cap \B^{<\omega} (\D)$ is pointwise closed.

Then setting 
\begin{equation*}
\zeta_{\D}^{\F} =  \xi_{\D}^{\F_h} = \sup 
\{ s_{\D_1} (\F_h) : \D_1 <\D\}\ ,
\end{equation*}
which is a countable ordinal, the following subcases obtain:
\begin{itemize}
\item[2(i)] 
If $\xi+1 <\zeta_{\D}^{\F}$, then there exists $\D_2<\D$ 
such that 
$$\B^\xi \cap \B^{<\omega} (\D_2) \subseteq \F\ ;$$ 
\item[2(ii)] 
if $\xi+1>\xi>\zeta_{\D}^{\F}$, then for every $\D_1<\D$ 
there exists $\D_2<\D_1$  such that 
$$\B^\xi\cap \B^{<\omega} (\D_2) \subseteq \B^{<\omega} (\nat)\setminus \F$$ 
(equivalently $\F\cap \B^{<\omega} (\D_2) \subseteq(\B^\xi)^*\setminus 
\B^\xi$); and 
\item[2(iii)] 
if $\xi+1 = \zeta_{\D}^{\F}$ or $\xi = \zeta_{\D}^{\F}$, then there 
exists $\D_2 <\D$ such that 
\item[{}] either ${\B_\xi\cap \B^{<\omega} (\D_2)\subseteq \F}$ or 
$\B_\xi\cap \B^{<\omega} (\D_2)\subseteq \B^{<\omega}(\nat)\setminus\F$.
\end{itemize}
\end{thm}

\begin{proof}{} 
[Case 1]
If the hereditary family $\F_h\cap \B^{<\omega} (\D)$ is not pointwise 
closed, then there exists $\D_1<\D$ such that $\B^{<\omega} (\D_1) 
\subseteq \F_h \subseteq \F$, according to 
Proposition~\ref{prop:finitefamily}.

\noindent [Case 2] 
If $\F_h \cap \B^{<\omega}(\D)$ is pointwise closed, then $\zeta_{\D}^{\F}$ 
is a countable ordinal, according to Theorem~\ref{block-NashWilliams}. 

2(i) 
Let $\xi+1<\zeta_{\D}^{\F}$. 
Then $\xi+1 <\xi_{\D}^{\F_h}$, and, according to 
Theorem~\ref{block-NashWilliams} (subcase 2(i)), there exists 
$\D_2 <\D$ such that 
\begin{equation*} 
\B^\xi \cap \B^{<\omega} (\D_2) \subseteq \F_h \subseteq \F\ .
\end{equation*}

2(ii)
Let $\xi +1>\xi >\zeta_{\D}^{\F}$ and $\D_1 <\D$. 
According to Theorem~\ref{block-NashWilliams} (subcase 2(ii)), there exists 
$\D'_1 <\D_1$ such that 
\begin{equation}\label{eq:block-NashWilliams} 
\B^{\zeta_{\D}^{\F}} \cap \B^{<\omega} (\D'_1) \subseteq \B^{<\omega} 
(\nat) \setminus \F_h\ .
\end{equation}
According to the $\xi$-block Ramsey partition theorem 
(Theorem~\ref{thm:block-Ramsey}),
there exists $\D_2 <\D'_1$ such that 
\begin{itemize}
\item[{}] either $\B^\xi \cap \B^{<\omega} (\D_2) \subseteq \F$ 
\quad or\quad  
$\B^\xi\cap \B^{<\omega} (\D_2) \subseteq \B^{<\omega} (\nat) 
\setminus\F$.
\end{itemize}
We claim that the first alternative does not hold. 
Indeed, if $\B^\xi \cap \B^{<\omega}(\D_2)\subseteq \F$, then 
$(\B^\xi \cap \B^{<\omega} (\D_2))^* \subseteq \F^* = \F$. 
Using the canonical representation of every infinite disjoint collection 
with respect to $\B_\xi$ (Proposition~\ref{prop:canonicalrep}) it is 
easy to check that 
\begin{equation*}
(\B^\xi)^* \cap \B^{<\omega} (\D_2) = (\B^\xi\cap \B^{<\omega} (\D_2))^*\ .
\end{equation*}
Hence, $(\B^\xi)^* \cap \B^{<\omega} (\D_2)\subseteq \F$.

Since $\xi>\zeta_{\D}^{\F}$, there exists $\D_3 <\D_2$ such that 
\begin{equation*}
(\B^{\zeta_{\D}^{\F}})_* \cap \B^{<\omega} (\D_3) \subseteq (\B^\xi)^* 
\cap \B^{<\omega} (\D_2) \subseteq \F\ ,
\end{equation*}
(Corollary~\ref{cor:tree}). 
So, $(\B^{\zeta_{\D}^{\F}})_* \cap \B^{<\omega} (\D_3) \subseteq \F_h$. 
This is a contradiction to the relation \eqref{eq:block-NashWilliams}. 
Hence, $\B^\xi \cap \B^{<\omega} (\D_2) \subseteq \B^{<\omega}(\nat) 
\setminus \F$ and equivalently $\F\cap \B^{<\omega} (\D_2)\subseteq 
(\B_\xi)^* \setminus \B_\xi$ (Proposition~\ref{prop:tree}). 

2(iii)
In the cases $\zeta_{\D}^{\F} = \xi+1$ or $\zeta_{\D}^{\F} = \xi$, we 
use Theorem~\ref{thm:block-Ramsey}.
\end{proof}

\begin{cor}\label{cor:tree2}
Let $\F$ be a family of finite disjoint collections on $\nat$ which is a 
tree and let $\D$ be an infinite disjoint collection on $\nat$. 
Then
\begin{itemize}
\item[(i)] either 
there exists $\D_1 <\D$ such that $\B^{<\omega} (\D_1)\subseteq \F$;
\item[(ii)] or there exists a countable ordinal $\xi_0 = \zeta_{\D}^{\F}$ 
such that for all $\xi >\xi_0$ and $\D_1 <\D$ there exists 
$\D_2 <\D$ such that $\B^\xi \cap \B^{<\omega}(\D_2) \subseteq \B^{<\omega} 
(\nat) \setminus \F$.
\end{itemize}
\end{cor}

\begin{proof} 
We apply Theorem~\ref{block-NashWilliams2} (case 1 and subcase 2(ii)).
\end{proof}

Condition 2(ii) of Theorems~\ref{block-NashWilliams} and 
\ref{block-NashWilliams2} (and (ii) of Corollary~\ref{cor:tree2}) has 
a number of interesting reformulations.

\begin{prop}\label{prop:canonicalrep2}
Let $\F\subseteq \B^{<\omega} (\nat)$, $\xi<\omega_1$ and 
$\D\in \B^{<\omega} (\nat)$.
The following are equivalent:
\begin{itemize}
\item[(i)] $\B^\xi \cap \B^{<\omega} (\D) \subseteq \F$.
\item[(ii)] For every $\D_1 <\D$ the unique initial segment of $\D_1$ 
which is an element of $\B^\xi$ belongs to $\F$.
\item[(iii)] Given any  sequence $(\D_n)_{n\in\nat}$ of infinite 
disjoint collections with $\D_n <\D$ for  every $n\in\nat$, and any 
$s_n \in FU(\D_n)$ for every $n\in\nat$ with $s_1<s_2<\cdots$,
%$s_2\in FU(\D_2)$ with $s_1<s_2$, $s_3\in FU(\D_3)$, 
then  there exists 
$n_0\in\nat$ such that $(s_1,\ldots,s_{n_0}) \in \B^\xi \cap \F$.
\end{itemize}
\end{prop}

\begin{proof} 
(i) $\Rightarrow$ (ii). 
Let $\D_1 <\D$. 
Using the canonical representation of $\D_1$ with respect to $\xi$
(Proposition~\ref{prop:canonicalrep}) there exists a unique initial 
segment $\bs$ of $\D_1$ which is an element of $\B^\xi$. 
Since $\bs \in \B^\xi \cap \B^{<\omega} (\D)$, we have that 
$\bs \in \F$, according to (i).
\medskip

\noindent (ii) $\Rightarrow$ (i). 
Let $\bs \in \B^\xi \cap \B^{<\omega} (\D)$, where $\D = (s_n)_{n\in\nat}$
with $s_1<s_2<\cdots$ and $\bs = (t_1,\ldots,t_k)$ with 
$t_1<\cdots < t_k\in FU(\D)$. 
Then $\D_1 = \{t_1,\ldots,t_k\} \cup (\D- t_k)$ is an 
infinite disjoint collection and of course $\D_1 <\D$. 
Using the canonical representation of $\D_1$ with respect to $\B_\xi$ 
(Proposition~\ref{prop:canonicalrep}), and (ii), we have that $\bs \in \F$.
\medskip

\noindent (i) $\Rightarrow$ (iii). 
Let $(\D_n)_{n\in\nat}$ be a sequence of 
infinite disjoint collections with $\D_n<\D$ for every $n\in\nat$. 
Choose $s_n \in FU(\D_n)$ for every $n\in\nat$ with $s_1<s_2<\cdots$.
The infinite disjoint collection $(s_n)_{n\in\nat}$ has canonical 
representation with respect to $\B_\xi$, hence there exists unique $n_0\in\nat$
such that $(s_1,\ldots,s_{n_0})\in \B^\xi \cap \B^{<\omega} (\D)$. 
According to (i), we have that $(s_1,\ldots,s_{n_0})\in \F$.
\medskip

\noindent (iii) $\Rightarrow$ (ii). 
Let $\D_1<\D$ and let $\D_1= (t_n)_{n\in\nat}$. 
Using (iii) with $\D_n = \D_1$ for every $n\in\nat$ and $s_n=t_n$ for 
every $n\in\nat$ we have the existence of an $n_0\in\nat$ such that 
${(t_1,\ldots,t_{n_0}) \in \B^\xi \cap \F}$.
\end{proof}

We express now Corollary~\ref{cor:tree2}, in terms of sumsets.

\begin{cor}\label{cor:sumsets}
Let $\F\subseteq [\nat]^{<\omega}$ be a family of finite subsets of $\nat$ 
which is a tree. 
Then
\begin{itemize}
\item[{}] either there exists $L\in[\nat]$ such that $P(L)\subseteq \F$;
\item[{}] or there exists a countable ordinal $\xi_0 = \xi(\F)$ such that 
for all $\xi>\xi_0$ there exists $L\in[\nat]$ such that $P^\xi(L)
\subseteq [\nat]^{<\omega}\setminus \F$.
\end{itemize}
\end{cor}

\begin{proof}
We apply Theorem~\ref{block-NashWilliams2} for the family 
$\F_1 = \bar\varphi^{-1} (\F)$, where $\bar\varphi:\B^{<\omega}(\nat) \to 
[\nat]^{<\omega}$ with $\bar\varphi ((s_1,\ldots,s_k)) = 
(\varphi (s_1),\ldots,\varphi (s_k))$. 
If $\F$ is a tree, then obviously $\F_1$ is a tree in $\B^{<\omega}(\nat)$.
\end{proof}

We now consider a very simplified statement of 
Theorem~\ref{block-NashWilliams2} 
(not including countable ordinals), more akin to the Gowers' 
reformulation (in \cite{G}) of the classical Nash-Williams theorem
(\cite{NW}).
%In fact this statement is equivalent if trees can be replaced by 
%arbitrary partitions of $\B^{<\omega}(\nat)$.

\begin{cor}%[Block Nash-Williams partition theorem in Gower's reformulation]
\label{cor:gowers}
Let $\F\subseteq \B^{<\omega}(\nat)$ be a family of finite disjoint 
collections %(which is a tree) 
and $\D$ an infinite disjoint collection.
Then there exists $\D_1<\D$ such that 
\begin{itemize}
\item[(i)] either $\B^{<\omega} (\D_1)\subseteq \F$; 
\item[(ii)] or for every $\D_2 <\D_1$ there exists an initial segment of $\D_2$
which belongs to \newline $\B^{<\omega} (\nat)\setminus \F$.
\end{itemize}
\end{cor}

\begin{proof}
If $\F$ is a tree, then the proof follows from 
Theorem~\ref{block-NashWilliams2}. 

Suppose $\F$ is not a tree. 
Then we set $\F_t = \{\bs \in \F:$ every $\bt\propto\bs$ belongs to 
$\F\}\cup \{\emptyset\}$. 
The family $\F_t$ is a tree and $\F_t \subseteq \F$. 
According to Theorem~3.3, there exists $\D_1 <\D$ such that 
\begin{itemize}
\item[{}] either $\B^{<\omega} (\D_1) \subseteq \F_t \subseteq \F$; 
\item[{}] or for every $\D_2<\D_1$ there exists an initial segment of $\D_2$ 
which belongs to \newline 
$\B^{<\omega} (\nat)\setminus \F_t$.
\end{itemize}

Let $\D_2 <\D_1$, and let $\bs \propto \D_2$ with $\bs \in\B^{<\omega}
(\nat)\setminus \F_t = (\F\setminus \F_t)\cup (\B^{<\omega}(\nat)
\setminus \F)$. 
Then, either $\bs \in \B^{<\omega}(\nat)\setminus \F$, as required, 
or $\bs\in \F\setminus \F_t$. 
In case $\bs\in \F\setminus \F_t$, by the definition of $\F_t$, there 
exists $\bs_1\propto \bs$ such that $\bs_1\in \B^{<\omega} (\nat) \setminus\F$,
as required.
\end{proof}

Corollary~\ref{cor:gowers} implies the following 
block Nash-Williams theorem, proved by Milliken (\cite{M}).
%(and is a special case of Carlson's general partition theorem \cite{C}). 
In fact, Corollary~3.7 is equivalent to this theorem.

\begin{cor} 
[Block Nash-Williams theorem, \cite{M}]
\label{cor:blockNW}
Let $\U\subseteq \B^{\omega}(\nat)$ be a family of infinite disjoint 
collections, which is closed in the pointwise topology of 
$\{0,1\}^{[\nat]^{<\omega}}$ and $\D$ an infinite disjoint collection. 
Then there exists $\D_1 <\D$ such that 
$$\text{either }\ \B^\omega (\D_1)\subseteq \U\ ;\quad \text{or}\quad 
\B^\omega (\D_1)\subseteq \B^\omega (\nat) \setminus \U\ .$$
\end{cor}

\begin{proof} 
Let $\F_{\U} = \{\bs \in \B^{<\omega}(\nat)$: there exists $\D'\in \U$ 
such that $\bs\propto \D'\}$.
Since the family $\F_{\U}$ is a tree, we use Corollary~\ref{cor:gowers}.
Then we have the following two cases:

\noindent {\bf [Case 1]}
There exists $\D_1<\D$ such that $\B^{<\omega}(\D_1)\subseteq \F_{\U}$. 
Then, $\B^\omega (\D_1)\subseteq \U$.
Indeed, if $\D_2 = (s_n)_{n\in\nat} \in \B^\omega(\D_1)$, then 
$(s_1,\ldots,s_n)\in \F_{\U}$ for every $n\in\nat$. 
Hence, for every $n\in\nat$ there exists $\D'_n\in\U$ such that 
$(s_1,\ldots,s_n)\in \D'_n$. 
Since $(\D'_n)_{n\in\nat}$ converges pointwise to $\D_2$ and $\U$ is closed, 
we have that $\D_2\in\U$ and consequently that $\B^\omega(\D_1)\subseteq\U$.

\noindent {\bf [Case 2]}
There exists $\D_1<\D$ such that for every $\D_2 <\D_1$ there exists an 
initial segment $(s_1,\ldots,s_n)$ of $\D_2$ which belongs to 
$\B^{<\omega} (\nat)\setminus \F_{\U}$. 
Then, $\D_2 \in \B^\omega (\nat)\setminus \U$ for every $\D_2<\D_1$. 
Hence, $\B^\omega (\D_1) \subseteq \B^\omega (\nat)\setminus \U$.
\end{proof}

\begin{remark}\label{rem:block-NashWilliams}
The block Nash-Williams partition theorem (Corollary~\ref{cor:blockNW})
is in fact equivalent to Corollary~\ref{cor:gowers}. 
In fact, let $\F\subseteq \B^{<\omega} (\nat)$. 
Set 
$$\U_{\F} = \{ \D\in \B^{<\omega} (\nat): \text{ there exists }
(s_1,\ldots,s_k)\in\F\text{ such that } (s_1,\ldots,s_k)\propto \D\}\ .$$
The complement $\B^\omega (\nat)\setminus \U_{\F}$ is pointwise closed, 
so, using Corollary~\ref{cor:blockNW} for the family 
$\B^\omega(\nat)\setminus \U_{\F}$, we obtain Corollary~\ref{cor:gowers}. 
On the other hand, Corollary~\ref{cor:gowers}, as we have already seen, 
implies block Nash-Williams (Corollary~3.8).

%\begin{prop}\label{prop3.10}
%{\tt [Cited on page  14]}
%\end{prop}

%%%%%%%%  4444444 %%%%%%%%%%%
\section{The Block Ellentuck partition theorem for every countable ordinal}

In this section we show that our main theorem~\ref{block-NashWilliams2} 
implies a block partition theorem, for 
partitions of all infinite block sequences closed in the analogue of the 
Ellentuck topology $\Tau_E$, and for every 
countable ordinal (Theorem~\ref{thm:Ellentuck}). 
A simple consequence of Theorem~\ref{thm:Ellentuck} (together with 
Corollary~\ref{cor:Ellentuck2}) is the characterization of 
completely Ramsey block
families in terms of the Baire property in the topology $\Tau_E$, a result 
proved with different methods by Milliken \cite{M}. 
%(and also by Carlson \cite{C}). 

We define below the topology $\Tau_E$ on $\B^\omega (\nat)$, a 
block analogue of the Ellentuck topology (\cite{E}).

\begin{defn}\label{def:Ellentuck}
Let $\Tau_E$ be a topology of $\B^\omega(\nat)$ with basic open sets of 
the following form:
\begin{equation*}
[\bs ,\D] =  \{\D_1 \in \B^\omega (\nat): \bs \propto \D_1\quad\text{and}
\quad \D_1 - \cup_{F\in \bs} F <\D\}\ ,
\end{equation*}
where $ \bs \in \B^{<\omega}(\nat)$ and
$\D\in \B^\omega (\nat)$.

The topology $\Tau_E$ is stronger than the relative topology of 
$\B^{<\omega}(\nat)$ with respect of the pointwise convergence topology of 
$\{0,1\}^{[\nat]^{<\omega}}$, which has basic open sets of the form 
$[\bs,\nat] = \{\D_1 \in \B^\omega (\nat): \bs\propto \D_1\}$.

We denote by $\hat\U$ and $\U^\lozenge$ the closure and the interior 
respectively of a family $\U\subseteq \B^\omega(\nat)$ in the topology 
$\Tau_E$. 
Then it is easy to see that 
\begin{equation*}
\begin{split}
&\hat\U = \{\D\in\B^\omega (\nat): [\bs ,\D]\cap \U\ne\emptyset
\ \text{ for every }\  \bs\propto \D\}\ ;\ \text{ and}\\
&\U^\lozenge = \{\D\in\B^\omega (\nat):\ \text{ there exists }
\bs\propto \D\ \text{ such that }\ 
[\bs,\D] \subseteq \U\}\ .
\end{split}
\end{equation*}
\end{defn}
\end{remark}

\begin{thm}\label{thm:Ellentuck}
Let $\U\subseteq \B^\omega (\nat)$ be a family of infinite disjoint 
collections, $\bs\in\B^{<\omega}(\nat)$ and $\D$ an infinite disjoint 
collection. 
Then
\begin{itemize}
\item[{}] either there exists $\D_1<\D$ such that $[\bs,\D_1]\subseteq \hat\U$;
\item[{}] or there exists a countable ordinal $\xi_0= \zeta_{(\bs,\D)}^{\U}$ 
such that for every $\xi>\xi_0$ there exists $\D_1<\D$ with 
$[\bs \cup \bt,\D_1]\subseteq \B^\omega (\nat)\setminus \U$ for every 
$\bt \in \B^\xi \cap \B^{<\omega} (\D_1)$.
\end{itemize}
\end{thm}

We will give the proof of this theorem after the following lemma which 
is analogous to Lemma~\ref{lem:block-Ramsey}. 

\begin{lem}\label{lem:Ellentuck}
Let $\G \subseteq\{[\bs,\D] :\bs\in \B^{<\omega}(\nat)$ and 
$\D\in \B^\omega (\nat)\}$ with the following two properties:
\begin{itemize}
\item[(i)] for every $(\bs,\D)\in \B^{<\omega}(\nat)\times\B^\omega (\nat)$
there exists $\D_1<\D$ such that $[\bs,\D_1]\in\G$; and 
\item[(ii)] for every $[\bs,\D_1]\in\G$ and $\D_2<\D_1$ we have 
$[\bs,\D_2] \in\G$.
\end{itemize}
Then, for every $(\bs,\D)\in \B^{<\omega} (\nat)\times \B^\omega(\nat)$ 
there exists $\D_0\in [\bs,\D]$ such that $[\bs\cup\bt,\D_1]\in\G$ 
for every $\bt\in \B^{<\omega}(\D_0)$ with 
$\bigcup_{F\in\bs} F < \bigcup_{F\in\bt} F$ if $\bt \ne\emptyset$ 
and $\D_1<\D_0$.
\end{lem}

\begin{proof}
Let $(\bs,\D)\in\B^{<\omega}(\nat)\times \B^\omega(\nat)$.
We can assume that $\bigcup_{F\in\bs} F< \bigcup_{F\in\D} F$. 
According to assumption~(i), there exists $\D_1<\D$ such that 
$[\bs,\D_1] \in\G$. 
Let $\D_n<\cdots <\D_1 \in\B^\omega (\nat)$ have been constructed and 
$\D_n = (s_i^n)_{i\in\nat}$ for every $n\in\nat$.

Set $\{\bt_1,\ldots,\bt_r\} = \{\bt\in\B^{<\omega}(\nat) :\bt\subseteq 
FU (\{s_1^1,\ldots, s_n^n\})\}$. 
According to (i), there exist $\D_{n+1}^1 <\D_n - s_n^n$ such that 
$[\bs\cup \bt_1,\D_{n+1}^1]\in\G$, $\D_{n+1}^2 <\D_{n+1}^1$ such that 
$[\bs\cup\bt_2,\D_{n+1}^2]\in \G$, and finally $\D_{n+1}^r < \D_n-s_n^n$ 
such that $[\bs\cup \bt_r,\D_{n+1}^r]\in \G$. 
Set $\D_{n+1} = \D_{n+1}^r = (s_i^{n+1})_{i\in\nat}$. 
Then, according to (ii), $[\bs\cup \bt_i,\D_{n+1}]\in\G$ for every 
$1\le i\le r$. 

Set $\D_0 = \{s_n^n :n\in\nat\} \cup \bs$. 
Then $\D_0 \in [\bs,\D]$.
Let $\bt \in \B^{<\omega} (\D_0)$ with $t\ne\emptyset$ and 
${\bigcup_{F\in\bs} F < \bigcup_{F\in\bt} F}$. 
If $n_0 = \min \{n\in\nat: \bt \subseteq FU(\{s_1^1,\ldots,s_n^n\})\}$,
then 
$[\bs\cup\bt,\D_{n_0+1}]\in \G$. 
According to assumption~(ii), $[\bs\cup \bt,\D_0 - s_{n_0}^{n_0}]\in \G$. 
Hence, $[\bs\cup \bt,\D_0] = [\bs\cup\bt,\D_0 - s_{n_0}^{n_0}]\in \G$.
If $\bt =\emptyset$, then $[\bs,\D_1]\in\G$, hence $[\bs,\D_0] \in \G$.
\end{proof}

\begin{proof}[Proof of Theorem~\ref{thm:Ellentuck}] 
Let $\U\in \B^\omega (\nat)$, $\bs\in \B^{<\omega}(\nat)$ and 
$\D \in \B^\omega (\nat)$. 
Set 
\begin{equation*}
\begin{split}
\G_{\U} = &\{[\bs,\D] :(\bs,\D) \in\B^{<\omega}(\nat)\times\B^\omega (\nat)
\ \text{ and}\\
&\qquad \text{either }\ [\bs,\D]\cap U = \emptyset\quad\text{or}\quad 
[\bs,\D_1]\cap U\ne\emptyset\ \text{ for every }\ \D_1<\D\}\ .
\end{split}
\end{equation*}
It is easy to check that $\G_U$ satisfies the assumptions~(i) and (ii) 
of Lemma~\ref{lem:Ellentuck}, 
hence, there exists $\D_0\in[\bs,\D]$ such that $[\bs\cup \bt,\D_0]\in\G_{\U}$
for every $\bt \in \B^{<\omega} (\D_0)$ with 
$\bigcup_{F\in\bs}F < \bigcup_{F\in\bt} F$ if $\bt \ne\emptyset$.
Set
\begin{equation*}
\F = \{\bt\in\B^{<\omega} (\D_0) : \cup_{F\in\bs}F < \cup_{F\in\bt} F
\text{ if } \bt\ne\emptyset\text{ and } 
[\bs \cup \bt ,\D_1] \cap U\ne \emptyset\text{ for every } \D_1<\D_0\}\ .
\end{equation*}
The family $\F$ is a tree. 
Indeed, let $\bt \in \F$ and $\bt_1 \propto\bt$.
Then $[\bs\cup \bt_1,\D_0] \in\G_{\U}$.
So either $[\bs \cup \bt_1,\D_0]\cup \U=\emptyset$ which is impossible, since 
$[\bs \cup \bt,\D_0]\cap U\ne\emptyset$, or $[\bs\cup \bt_1,\D_1]\cap U
\ne\emptyset$ for every $\D_1<\D_0$. 
Hence $\bt_1 \in \F$. 

We use Theorem~\ref{block-NashWilliams2}, for $\F$ and 
$\D_0 - \cup_{F\in \bs}F$.
We have the following cases:

\noindent {\bf [Case 1]} 
There exists $\D_1<\D_0 - \bigcup_{F\in\bs} F <\D$ such that 
$\B^{<\omega} (\D_1) \subseteq \F$. 
This gives that $[\bs\cup \bt,\D_2]\cap \U\ne \emptyset$ for every 
$\bt \in \B^{<\omega} (\D_1)$ and $\D_2 <\D_1$, 
which implies that $[\bs,\D_1] \subseteq\hat\U$.

\noindent {\bf [Case 2]} 
There exists a countable ordinal $\xi_0 = \zeta_{\D}^{\F} 
= \zeta_{(\bs,\D)}^{\U}$ such that for every $\xi>\xi_0$ there exists 
$\D_1 <\D_0 - \bigcup_{F\in \bs} F <\D$ with 
$\B^\xi \cap \B^{<\omega} (\D_1) \subseteq \B^{<\omega}(\nat)\setminus\F$.
Then $[\bs\cup \bt,\D_1] \subseteq \B^\omega (\nat)\setminus \U$ for every 
$\bt \in \B^\xi \cap \B^{<\omega} (\D_1)$.
\end{proof}

Applying Theorem~\ref{thm:Ellentuck} to partitions $\U$ that are closed 
(the pointwise closed families are included in this class) or meager in the 
topology $\Tau_E$, we consider the following consequences.

\begin{cor}\label{cor:Ellentuck}
Let $\U$ be a subset of $\B^\omega (\nat)$ pointwise closed (or, even $\U$ 
closed in the topology $\Tau_E$), $\bs \in \B^{<\omega} (\nat)$ and $\D$ 
an infinite disjoint collection. 
Then 
\begin{itemize}
\item[{}] either there exists $\D_1<\D$ such that $[\bs,\D_1]\subseteq \U$; 
\item[{}] or there exists a countable ordinal $\xi_0 = \xi_{(\bs,\D)}^{\U}$, 
such that for every $\xi>\xi_0$ there exists $\D_1<\D$ such that 
$[\bs \cup \bt,\D_1] \subseteq\B^\omega (\nat)\setminus \U$ for every
$\bt\in\B^\xi \cap \B^{<\omega} (\D_1)$.
\end{itemize}
\end{cor}

\begin{proof} 
It follows from Theorem~\ref{thm:Ellentuck}.
\end{proof}

\begin{cor}\label{cor:Ellentuck2}
Let $\U$ be a subset  of $\B^\omega (\nat)$ meager in the topology $\Tau_E$, 
$\bs\in \B^{<\omega} (\nat)$ and $\D\in\B^\omega (\nat)$. 
Then, there exists a countable ordinal $\xi_0$ such that for every $\xi>\xi_0$
there exists $\D_1 <\D$ such that $[\bs\cup \bt,\D_1] \subseteq \B^\omega
(\nat)\setminus \U$ for every $\bt\in \B^\xi \cap \B^{<\omega} (\D_1)$.
\end{cor}

\begin{proof}
We use Theorem~\ref{thm:Ellentuck} for $\U$.
We will prove that the first alternative is impossible. 
Indeed, let $\D_1<\D$ such that $[\bs,\D_1] \subseteq\hat\U$. 
If $\U= \bigcup_{n\in\nat} \U_n$ with $(\hat\U_n)^\lozenge = \emptyset$ 
for every $n\in\nat$ and $|\bt|$ is the cardinality of $\bt\in\B^{<\omega)} 
(\nat)$, then set 
\begin{equation*}
\begin{split}
\G = &\{[\bt,\D_0] : \bt \in \B^{<\omega} (\nat),\ \D_0\in \B^\omega(\nat)
\ \text{ and}\\
&\qquad [\bt,\D_0]\cap \U_k =\emptyset\ \text{ for every } \ k\in\nat
\text{ with } k\le |\bt|\}\ .
\end{split}
\end{equation*}
The family $\G$ satisfies the conditions (i) and (ii) of 
Lemma~\ref{lem:Ellentuck}. 
Indeed,  $\G$ satisfies (i), since,
according to Theorem~\ref{thm:Ellentuck}, for every $\bt\in\B^{<\omega}
(\nat)$, $\D\in\B^\omega (\nat)$ and $k\in\nat$ there exists 
$\D'_0<\D$ such that $[\bt,\D'_0]\cap \U_k = \emptyset$, as 
it is impossible $[\bt,\D'_0]\subseteq \hat\U_k$  and 
obviously $\G$ satisfies (ii). 
Hence, according to Lemma~\ref{lem:Ellentuck}, there exists 
$\D_2 \in [\bs,\D_1]$ such that 
$[\bs\cup \bt,\D_2]\in\G$ for every $\bt\in \B^{<\omega}(\D_2)$. 
Then, $[\bs,\D_2]\cap \U=\emptyset$. 
Indeed, let $\D_3\in [\bs,\D_2]\cap\U$. 
Then, $\D_3 \in [\bs,\D_2] \cap \U_k$ for some $k\in\nat$. 
Hence, there exists $\bt\in\B^{<\omega} (\D_2)$ with $\bs\cup \bt\propto\D_3$,
$k\le |\bs\cup \bt|$ and $[\bs\cup\bt,\D_2]\cap \U_k\ne\emptyset$.
Then, $[\bs\cup \bt,\D_2]\notin \G$. 
A contradiction, since $\bt\in\B^{<\omega}(\D_2)$. 
Thus, $[\bs,\D_2]\cap \U=\emptyset$.
This implies that $\D_2\notin \hat\U$, which is a contradiction, since 
$\D_2 \in [\bs ,\D_1] \subseteq\hat\U$. 
Hence, the first alternative of Theorem~\ref{thm:Ellentuck} for $\U$ is 
impossible, so the second alternative holds for $\U$.
\end{proof}

We recall the definition of the completely Ramsey families of infinite
disjoint collection, given in \cite{M}.

\begin{defn}\label{def:Milliken}
A family $\U$ of infinite disjoint collections is called {\em completely 
Ramsey} if for every $\bs\in\B^{<\omega} (\nat)$ and every infinite 
disjoint collection $\D$, there exists $\D_1<\D$ such that 
\begin{equation*}
\text{either }\ [\bs,\D_1]\subseteq \U\quad\text{ or }\quad 
[\bs,\D_1]\subseteq \B^\omega (\nat) \setminus \U\ .
\end{equation*}
\end{defn}

A further consequence of Theorem~\ref{thm:Ellentuck} is the following 
corollary which gives the characterization of completely Ramsey families 
(proved also in \cite{M}).

\begin{cor}\label{cor:Baireproperty}
A family $\U\subseteq \B^\omega (\nat)$ is completely Ramsey if and only if 
$\U$ has the Baire property in the topology $\Tau_E$.
\end{cor} 

\begin{proof} 
Let $\U\subseteq \B^\omega (\nat)$ has the Baire property in the topology 
$\Tau_E$. 
Then $\U = \B \triangle  \C  = (\B \cup \C^c) \cup (\C\cap \B^c)$, where 
$\B\subseteq \B^\omega (\nat)$ is $\Tau_E$-closed and 
$\C\subseteq\B^\omega (\nat)$ is $\Tau_E$-meager and 
$\C^c= \B^\omega (\nat)\setminus\C$.  
According to Corollaries~\ref{cor:Ellentuck} and \ref{cor:Ellentuck2}, 
for every $\bs\in \B^{<\omega}(\nat)$ and $\D\in\B^\omega (\nat)$, 
there exists $\D_1<\D$ such that $[\bs,\D_1]\subseteq \C^c$ and 
consequently there exists $\D_2<\D_1$ such that 
\begin{itemize}
\item[{}] either $[\bs,\D_2] \subseteq \B\cap [\bs,\D_1]\subseteq 
\B\cap \C^c \subseteq \U$
\item[{}] or $[\bs,\D_2]\subseteq \B^c \cap [\bs,\D_1] \subseteq 
\B^c \cap \C^c \subseteq \U^c$.
\end{itemize}
Hence, $\U$ is completely Ramsey.

On the other hand, if $\U$ is completely Ramsey, then $\U = \U^\lozenge
\cup (\U\setminus \U^\lozenge)$ and $\U\setminus\U^\lozenge$ is a meager set 
in $\Tau_E$. 
Hence $\U$ has the Baire property in the topology $\Tau_E$.
\end{proof}

\begin{remark}\label{rem:Ellentuck}
(i) Every subset $\U$ of $\B^\omega(\nat)$ which is  Borel set in the 
topology of pointwise convergence is completely Ramsey, since it has the 
Baire property in the topology $\Tau_E$.

(ii) Every subset of $\B^\omega(\nat)$ which is analytic set in the 
topology of pointwise convergence is completely Ramsey, since every 
analytic set (in this topology) has the Baire property.
\end{remark}

\section*{Acknowledgments}

The authors wish to thank Ted Odell and Haskell Rosenthal 
for helpful comments and discussions on the content of this paper. 
The authors wish to thank the Department of Mathematics of the University 
of Texas at Austin for the hospitality extended during their visit on the 
spring semester of 2004.

\bigskip
{\footnotesize
\noindent 
{\sc Department of Mathematics, Athens University, Athens 157 84, Greece}

\noindent 
{\sc Department of Mathematics, The University of Texas at Austin,
Austin, TX 78712 U.S.A.}

\noindent 
E-mail address of first named author: 

vfarmaki@math.uoa.gr, vfarmaki@math.utexas.edu

\noindent 
E-mail address of second named author: 

snegrep@math.uoa.gr, snegrep@math.utexas.edu

}

\end{document}